\documentclass[11pt]{amsart}
\usepackage{amssymb, paralist, xspace, graphicx, url, amscd, euscript, mathrsfs, stmaryrd}
\usepackage[all]{xy}
\SelectTips{cm}{}

\numberwithin{equation}{section}

1

\numberwithin{subsection}{section}

\allowdisplaybreaks[1]

\newenvironment{enumeratea}
{\begin{enumerate}[\upshape (a)]}
{\end{enumerate}}

\newenvironment{enumeratei}
{\begin{enumerate}[\upshape (i)]}
{\end{enumerate}}

\newenvironment{enumerate1}
{\begin{enumerate}[\upshape (1)]}
{\end{enumerate}}

\newtheorem*{namedtheorem}{\theoremname}
\newcommand{\theoremname}{testing}

\newtheorem{theorem}[subsection]{Theorem}
\newtheorem{proposition}[subsection]{Proposition}
\newtheorem{proposition-definition}[subsection]
{Proposition-Definition}
\newtheorem{corollary}[subsection]{Corollary}
\newtheorem{lemma}[subsection]{Lemma}

\theoremstyle{definition}
\newtheorem{definition}[subsection]{Definition}

\newtheorem{example}[subsection]{Example}

\newtheorem{remark}[subsection]{Remark}

\theoremstyle{remark}

\newcommand\nome{testing}
\newcommand\call[1]{\label{#1}\renewcommand\nome{#1}}
\newcommand\itemref[1]{\item\label{\nome;#1}}
\newcommand\refall[2]{\ref{#1}~(\ref{#1;#2})}
\newcommand\refpart[2]{(\ref{#1;#2})}

\newcommand\cB{\mathcal{B}}

\newcommand\cG{\mathcal{G}}

\newcommand\cI{\mathcal{I}}

\newcommand\cM{\mathcal{M}}
\newcommand\cN{\mathcal{N}}
\newcommand\cO{\mathcal{O}}
\newcommand\cP{\mathcal{P}}

\newcommand\cX{\mathcal{X}}
\newcommand\cY{\mathcal{Y}}

\newcommand\GG{\mathbb{G}}

\newcommand\ZZ{\mathbb{Z}}

\newcommand\bH{\mathbf{H}}

\newcommand\bL{\mathbf{L}}

\newcommand\rC{\mathrm{C}}

\newcommand\rZ{\mathrm{Z}}

\newcommand\rmm{\mathrm{m}}

\newcommand\frg{\mathfrak{g}}

\renewcommand\frm{\mathfrak{m}}

\newcommand\arr{\ifinner\to\else\longrightarrow\fi}

\newcommand\noqed{\renewcommand\qed{}}

\renewcommand\H{\operatorname{H}}

\newcommand\eqdef{\overset{\mathrm{\scriptscriptstyle def}} =}

\newcommand\into{\hookrightarrow}

\def\displaytimes_#1{\mathrel{\mathop{\times}\limits_{#1}}}

\def\displayotimes_#1{\mathrel{\mathop{\bigotimes}\limits_{#1}}}

\newcommand\ext{\operatorname{Ext}}

\newcommand\isom{\operatorname{Isom}}

\newcommand\aut{\operatorname{Aut}}

\newcommand\spec{\operatorname{Spec}}

\newcommand\id{\mathrm{id}}

\newcommand\pr{\operatorname{pr}}

\newdir{ >}{{}*!/-5pt/@{>}}

\newcommand\double{\rightrightarrows}

\newcommand\doublelong[2]{\mathbin{\xymatrix{{}\ar@<3pt>[r]^{#1}
\ar@<-3pt>[r]_{#2}&}}}

\newcommand{\underhom}
{\mathop{\underline{\mathrm{Hom}}}\nolimits}

\newcommand{\underisom}
{\mathop{\underline{\mathrm{Isom}}}\nolimits}

\newcommand{\underaut}
{\mathop{\underline{\mathrm{Aut}}}\nolimits}

\newlength{\ignora}
\newcommand{\hsmash}[1]{\settowidth{\ignora}{#1}#1\hspace{-\ignora}}

\newcommand{\qcoh}{\operatorname{QCoh}}

\newcommand{\coh}{\operatorname{Coh}}

\newcommand{\ind}{\operatorname{Ind}}

\newcommand{\catsch}[1]{(\mathrm{Sch}/#1)}

\newcommand{\cat}[1]{(\mathrm{#1})}

\newcommand{\qc}{quasi-coherent\xspace}

\newcommand{\lr}{linearly reductive\xspace}

\newcommand{\opp}{^{\mathrm{op}}}

\newcommand{\mmu}{\boldsymbol{\mu}}

\renewcommand{\aa}{\boldsymbol{\alpha}}

\newcommand{\red}{_{\mathrm{red}}}

\newcommand{\et}{_{\text{\'et}}}
\newcommand{\fppf}{_{\mathrm{fppf}}}

\newcommand{\homrep}{\underhom^{\mathrm{rep}}}
\newcommand{\hominj}{\underhom^{\mathrm{inj}}}

\newcommand{\thickslash}{\mathbin{\!\!\pmb{\fatslash}}}

\theoremstyle{plain}

\newtheorem{lem}[subsection]{Lemma}

\newcommand{\mc}{\mathcal}

\theoremstyle{definition}

\newcommand{\mls}{\mc }

\numberwithin{equation}{subsection}

\begin{document}

\title{Tame stacks in positive characteristic}

\author[Abramovich]{Dan Abramovich}

\author[Olsson]{Martin Olsson}

\author[Vistoli]{Angelo Vistoli}

\address[Abramovich]{Department of Mathematics\\
Brown University\\
Box 1917\\
Providence, RI 02912\\
U.S.A.}
\email{abrmovic@math.brown.edu}

\address[Olsson]{Department of Mathematics \#3840\\
University of California\\
Berkeley, CA 94720-3840\\
U.S.A.}
\email{molsson@math.berkeley.edu}

\address[Vistoli]{Scuola Normale Superiore\\Piazza dei Cavalieri 7\\
56126 Pisa\\ Italy}
\email{angelo.vistoli@sns.it}

\thanks{Vistoli supported in part by the PRIN Project ``Geometria
sulle variet\`a algebriche'', financed by MIUR. Olsson partially supported by NSF grant DMS-0555827 and an Alfred P. Sloan fellowship. Abramovich support in part by NSF grants DMS-0301695 and DMS-0603284.}

\date{March 11, 2007}

\maketitle


\section{Introduction}


Since their introduction in \cite{D-M, Artin}, algebraic stacks have been a key tool in the algebraic theory of moduli.  In characteristic 0, one often is able to work with Deligne--Mumford stacks, which, especially  in characteristic 0, enjoy a number of nice properties making them almost as easy to handle as algebraic spaces. In particular, if $\cM$ is a Deligne--Mumford stack and $M$ its coarse moduli space, then \'etale locally over $M$ we can present $\cM$ as the quotient of a scheme by a finite group action. In characteristic 0 the formation of $M$ commutes with arbitrary base change, and if $\cM \to S$ is flat then $M \to S$ is flat as well. A key property in characteristic 0 is that the pushforward functor $\qcoh(\cM) \to \qcoh(M)$ is exact.

In characteristics $p>0$ the situation is not as simple. First, in many situations one  needs to consider algebraic stacks (in the sense of Artin) with finite but possibly ramified diagonal. Examples include K3 surfaces, surfaces of general type, polarized torsors for abelian varieties, stable maps. Second, even Deligne--Mumford stacks may fail to have some desired properties, such as flatness of moduli spaces, as soon as the orders of stabilizers are divisible by the characteristics.

In this paper we isolate a class of algebraic stacks in positive and mixed characteristics, called \emph{tame stacks}, which is a good analogue to the class of Deligne--Mumford stacks in characteristic 0, and arguably better than the class of Deligne--Mumford stacks in positive and mixed characteristics. Their defining property is precisely the key property described above: an algebraic stack with finite diagonal is tame if and only if  the pushforward functor $\qcoh(\cM) \to \qcoh(M)$ is exact (see Definition \ref{Def:tame} and discussion therein for the precise hypotheses).

Our main theorem on tame stacks is Theorem \ref{thm:main}. In particular, we show that an algebraic stack $\cM$ is tame  if and only if, locally in the \'etale topology of the coarse moduli space $M$, it can be presented as the quotient $[U/G]$ of a scheme by the action of a finite flat \lr group scheme. Other desirable properties, such as flatness of coarse moduli spaces, commutation of the formation of moduli space with arbitrary base change, and stability of the tame property under pullbacks follow as corollaries.

It should be noted that it is significantly easier to show such a presentation as a quotient locally in the \emph{fppf} topology. But the presentation in the \'etale topology is extremely useful for the applications we envision. This is one of the most intricate technical points in the paper.

The proofs of our results on tame stacks necessitate a good classification and study of 
finite flat \lr group schemes. To our surprise, we have not found a good reference on this subject, and therefore developed it here. Our main theorem on \lr group schemes is Theorem \ref{Th:main-lr}. In particular it says that a finite flat group scheme $G\to S$ is \lr if and only if its geometric fibers $G_x$ are extensions $1 \to \Delta \to G_x \to Q \to 1$ with tame and \'etale quotient $Q$ and diagonalizable kernel $\Delta$. While this classification is pleasing in its simplicity, it is also a bit disappointing, as it shows that the only finite groups admitting a \lr reduction to characteristic $p>0$ have a normal and abelian $p$-Sylow subgroup.

Finite flat linearly reductive group schemes and their classification are the subject of Section~\ref{section1}.

In Section~\ref{section3} we define the notion of tame stack, and prove the key local structure theorem (Theorem~\ref{thm:main}).  

In a sequel to this paper, we develop the theory of twisted curves and twisted stable maps with a tame target, in analogy to \cite{A-V}. This is for us the main motivation for the introduction of tame stacks.  

In Appendix \ref{sec:rigidification} we discuss rigidification of stacks.  Discussion of rigidification has appeared in several places in the literature already, but unfortunately not in sufficient generality for the applications we have in mind.  So we treat the most general case here.

\subsection*{Acknowledgments} We would like to thank Andrew Kresch, Frans Oort and Ren\'e Schoof for useful discussions.

\section{Linearly reductive finite group schemes}\label{section1}

Throughout the paper all schemes are assumed to be quasi-separated. (Recall that a scheme $S$ is \emph{quasi-separated} when the diagonal $S \arr S \times S$ is quasi-compact.)

In this section \ref{section1} all group schemes will be flat, finite and finitely presented over an arbitrary scheme. Such a group scheme $G \arr S$ will be called \emph{constant} if $G$ is the product of $S$ by a finite group.

\subsection{Equivariant sheaves}

Let $\pi\colon G \arr S$ be a group scheme. A \emph{$G$-equivariant sheaf}  on $S$ 
 is a sheaf with an action of $G$. There are (at least) four ways of defining an action of $G$ on a \qc sheaf.

\begin{enumeratea}\call{definitions}

\itemref{1} A \qc sheaf $F$ on $S$ extends naturally to a functor
   \[
   F\colon \catsch{S}\opp \arr \cat{Grp}. 
   \]
If $f\colon T \arr S$ is a morphism of schemes, then we define $F(T)$ as $(f^{*}F)(T)$; each $F(T)$ has the structure of an $\cO(T)$-module.

Then an action of $G$ on $F$ is an $\cO$-linear action of the functor
   \[
   G \colon \catsch{S}\opp \arr \cat{Grp}
   \]
on $F$. In other words, for each $T \arr S$ we have an action of the group $G(T)$ on the $\cO(T)$-module $F(T)$, and this action is functorial in $T \arr S$.

\itemref{3} Same as above, but the action of $G(T)$ on $F(T)$ is only defined for flat locally finitely presented morphisms $T \arr S$.


\itemref{2} We have a sheaf of commutative Hopf algebras $\pi_{*} \cO_{G}$ on $S$. Then an action of $G$ on a \qc sheaf $F$ is defined as a coaction $F \arr F\otimes_{\cO_{S}}\pi_{*}\cO_{G}$ of this sheaf on $F$. Equivalently, in terms of the dual Hopf algebra $\bH_G = (\pi_*\cO_G)^\vee$, the ``convolution hyperalgebra of $G$", it is an action $F \otimes \bH_G \to F$.

\itemref{4} An action of $G$ on $F$ can also be defined as a isomorphism of $\cO_{G}$-modules $\phi\colon \pi^{*}F \simeq \pi^{*}F$ satisfying the cocycle  condition
   \[
   m^{*}\phi = \pr_{1}^{*}\phi \circ \pr_{2}^{*}\phi \colon
      \Pi^{*}F \simeq \Pi^{*}F,
   \]
where $m\colon G \times_{S} G \arr G$ is the multiplication morphism and $\Pi\colon G\times_{S}G \arr S$ is the structure morphism.
\end{enumeratea}

The connections among these various definitions are as follows. The equivalence between \refpart{definitions}{1} and \refpart{definitions}{4} is a particular case of \cite[Proposition~3.49]{descent}. The same Proposition, applied to the fibered category whose object are pairs $(T \to S, F)$, where $T \arr S$ is a flat finitely presented morphism, gives the equivalence between \refpart{definitions}{2} and \refpart{definitions}{4}. Finally, assigning a homomorphism of $\cO_{G}$-modules $\phi\colon \pi^{*}F \simeq \pi^{*}F$ is equivalent to assigning a homomorphism of $\cO_{S}$-modules $F \arr \pi_{*}\pi^{*}F = F \otimes_{\cO_{S}}\pi_{*}\cO_{G}$; and using \cite[Proposition~3.48 and Proposition~3.49]{descent} it is easy to see that $\phi$ satisfies the cocycle condition of \refpart{definitions}{4} if and only the corresponding homomorphism $F \arr F\otimes_{\cO_{S}}\pi_{*}\cO_{G}$ is a coaction.

Using any of the definitions above of a $G$-action, there is an obvious definition of an equivariant homomorphism of \qc sheaves with a $G$-action. We will denote by $\qcoh(S)$ the category of \qc sheaves on $S$, and by $\qcoh^{G}(S)$ the category of $G$-equivariant \qc sheaves over $S$.

When $S$ is locally noetherian, we also denote by $\coh(S)$ and $\coh^{G}(S)$ the categories of coherent sheaves, respectively without and with a $G$-action.

The category $\qcoh^{G}(S)$ is also equivalent to the category $\qcoh(\cB_{S}G)$ of \qc sheaves over the classifying stack $\cB_{S}G$. We will define this in a way that is slightly different from, but equivalent to, the standard one (\cite[Definition~13.2.2]{L-MB}). A \qc sheaf $F$ on $\cB_{S}G$ of $\qcoh^{G}(S)$ associates with each $G$-torsor $P \arr T$ an $\cO(T)$-module $F(P \to T)$; also, for each commutative diagram
   \begin{equation}\label{eq:torsors}
   \xymatrix{
   P' \ar[r]^{g}\ar[d] & P \ar[d]\\
   T' \ar[r]^{f}       & T
   }
   \end{equation}
where  the columns are $G$-torsors and $g$ is $G$-equivariant, we have a 
homomorphism $F(P \to T) \arr F(P' \arr T')$ that is linear with respect to the 
natural ring homomorphism $\cO(T) \arr \cO(T')$. These data  are required to 
satisfy the following conditions.

\begin{enumeratei}

\item Suppose that we are given a $G$-torsor $P \arr T$. Then we get a presheaf 
of $\cO_{S}$-modules $F_{P \to T}$ defined by sending a Zariski-open subscheme 
$U$ to the $\cO(U)$-module $F(P\mid_{U} \to U)$. We require  this to be a \qc 
sheaf on $S$.

\item Suppose that we have a commutative diagram like (\ref{eq:torsors}). Then 
we get a homomorphism of \qc sheaves $$F_{P\to T} \arr f_{*}F_{P' \to T'},$$ 
defined by the given homomorphism
   \[
   \xymatrix@R-20pt{
   F_{P \to T}(U)\ar@{=}[d] & F_{P' \to T'}\bigl(f^{-1}(U)\bigr) \ar@{=}[d]\\
   F(P\mid_{U} \arr U) \ar[r] & F\bigl (P'|_{f^{-1}(U)}\arr f^{-1}(U)\bigr)
   }
   \]
for each open subscheme $U \subseteq T$. Then the corresponding homomorphism 
$f^{*}{F_{P \to T}} \arr F_{P' \to T'}$ is required to be an isomorphism.
\end{enumeratei}

There is an obvious notion of homomorphism of \qc sheaves on $\cB_{S}G$: such a homomorphism $\phi\colon F \arr F'$ assigns to every $G$-torsor $P \to T$ a homomorphism of $\cO(T)$-modules $\phi_{P\to T}F(P \to T) \arr F'(P \to T)$, in such a way that given any commutative diagram (\ref{eq:torsors}), the diagram
   \[
   \xymatrix{
   F(P\to T)\ar[r]\ar[d]^{\phi_{P \to T}}
      & F(P'\to T')\ar[d]^{\phi_{P' \to T'}}\\
   F'(P\to T)\ar[r]
      & F'(P'\to T')
   }
   \]
commutes. The category $\qcoh(\cB_{S}G)$ is the category of \qc sheaves on $\cB_{S}G$, the arrows being homomorphisms.

The equivalence between $\qcoh(\cB_{S}G)$ and $\qcoh^{G}(S)$ is as follows.

Suppose that $F$ is an object of $\qcoh^{G}(S)$. If $T$ is an $S$-scheme and $P \arr T$ is a $G$-torsor, and $h\colon P \arr S$ is the composite of the given morphisms $P \arr T$ and $T \arr S$, then the pullback $h^{*}F$ is a \qc sheaf over $P$ with a $G$-action. On the other hand, by descent theory we have an equivalence between the category of $G$-equivariant \qc sheaves on $P$ and the category of \qc sheaves on $T$  (\cite[Theorem~4.46]{descent}); we define $\Phi(F)_{P\arr T}$ to be a quasi-coherent sheaf on $T$ whose pullback to $P$ is isomorphic to $h^{*}F$ a $G$-equivariant sheaf. It is easy to see the function that sends $P \arr T$ into $\Phi(F)_{P\arr T}$ has a natural structure of a \qc sheaf on $\cB_{S}G$, and that a homomorphism $f\colon F \arr F'$ of $G$-equivariant \qc sheaves on $S$ induces a homomorphism $\phi(f)\colon \Phi F \arr \Phi F'$ of \qc sheaves on $\cB_{S}G$. This defines the functor $\Phi\colon \qcoh^{G}(S)\arr \qcoh(\cB_{S}G)$.

Let us define the inverse functor $\Psi\colon \qcoh(\cB_{S}G) \arr \qcoh^{G}(S)$. Given a \qc sheaf $F$ on $\cB_{S}G$, we define the \qc sheaf $\Psi F$ on $S$ to be the sheaf $F_{G \to S}$ associated with the trivial $G$-torsor $G \arr S$. For each morphism $f\colon T \arr S$, the action of $G(T)$ on $f^{*}\Psi F$ is defined as follows: an element of $G(T)$ induces an an automorphism of the $G$-torsor $G_{T} \arr T$, which in turns induces an automorphism of $F_{G_{T} \to T} \simeq f^{*}\Psi F$. It is easy to see that $\Psi$ extends naturally to a functor $\Psi\colon \qcoh(\cB_{S}G) \arr \qcoh^{G}(S)$.

It is easy to see that the composite $\Phi\Psi$ is isomorphic to $\id_{\qcoh^{G}(S)}$. It is slightly less trivial to show that $\Psi\Phi$ is isomorphic to $\id_{\qcoh(\cB_{S}G)}$. The point is the following. Given a \qc sheaf $F$ on $\cB_{S}G$ and a $G$-torsor $\rho\colon P \arr T$ on an $S$-scheme $f\colon T \arr S$, the pullback $\pr_{2}\colon P \times_{T} P \arr P$ of $P \arr T$ to $P$ has canonical section: this induces a cartesian diagram
   \[
   \xymatrix{
   P \times_{T} P \ar[r] \ar[d]^{\pr_{2}}& G \ar[d]^{\pi}\\
   P \ar[r]^{f\rho} & S\hsmash{.}
   }
   \]
Since the pullback $\rho^{*}F_{P \to T}$ is isomorphic to $F_{P \times_{T} P \to P}$, this diagram induces an isomorphism of $\rho^{*}F_{P \to T}$ with $\rho^{*}f^{*}F_{G \to S}$. This isomorphism is easily seen to be $G_{T}$-equivariant; hence it descends to an isomorphism $F_{P \to T} \simeq f^{*}F_{G \to S} = (\Psi\Phi)F$.

Suppose that $\phi\colon H \arr G$ is a homomorphism of group schemes, there are two natural additive functors, the \emph{restriction functor}
   \[
   \phi^{*}\colon \qcoh^{G}(S) \arr \qcoh^{H}(S)
   \]
and the \emph{induction functor}
   \[
   \phi_{*}\colon \qcoh^{H}(S) \arr \qcoh^{G}(S).
   \]
The first is evident. The second, $\phi_*$ can be defined using functorial actions, or using the convolution algebras: $\phi_*F = F\otimes_{\bH_H}\bH_G$. It is also useful to think about them as follows: $\phi$ induces a morphism of algebraic stacks  
   \[
   \Phi\colon \cB_{S} H \arr \cB_{S}G
   \]
defined as usual by sending a principal $H$-bundle $Q \to T$ to the principal $G$-bundle $(Q \times_T G_T) / H_T$, the quotient taken with $H_T$ acting via $(q,g) \mapsto (qh, h^{-1}g)$. Then $\phi^{*}$ is pullback of \qc sheaves along $\Phi$, while $\phi_{*}$ is pushforward along $\Phi$.

\begin{remark}\call{rmk:important-points}
A few important points about these functors:

\begin{enumerate1}
\itemref{1} The functor $\phi^*$ is always exact. Indeed, in terms of actions, $\phi^*F$ is the same sheaf $F$ but with the $G$ action replaced by the action of $H$ through $\phi$, and the action does not intervene in exactness.

\itemref{2} If $H$ is a subgroup scheme of $G$, then $\Phi$ is finite, and in particular affine; hence $\phi_{*}$ is exact. In this case we denote it by $\ind^{G}_{H}$.

\itemref{3} 
If we think of the structure morphism $\pi\colon G \arr S$ as a homomorphism to the trivial group scheme and $F$ is a $G$-equivariant \qc sheaf on $S$, then we denote $\pi_{*}F$ by $F^{G}$. This \qc sheaf $F^{G}$ is naturally embedded in $F$ (by the adjunction map $\pi ^*\pi _*F\rightarrow F$),
and is called the \emph{invariant subsheaf}.

The invariant subsheaf can also be defined directly from any of the various definition of an action of $G$ on a \qc sheaf. For example, if $\alpha\colon F \arr F\otimes_{\cO_{S}}\pi_{*}\cO_{G}$ is a coaction, then $F^{G}$ is the kernel of $\alpha-\iota$, where $\iota\colon F \arr F\otimes_{\cO_{S}}\pi_{*}\cO_{G}$ is the trivial coaction, given by $s \mapsto s \otimes 1$.

\itemref{4} Suppose $\phi:H \to G$ is surjective, with kernel a flat group scheme $K$. For $F\in \qcoh^H(S)$ we have $\phi_*F = F^K$ with the induced action of $G$.  On the other hand if $F\in \qcoh^G(S)$ then the adjunction morphism $F \to \phi_*\phi^*F$ is an isomorphism, since the action of $K$ on $\phi^*F$ is trivial. In other words, we have a canonical isomorphism $ \phi_*\circ\phi^*\simeq \operatorname{id}$.
\end{enumerate1}
\end{remark}

%

\subsection{Linearly reductive group schemes}
\begin{definition} \label{Def:lr}
A group scheme $G \arr S$ is \emph{\lr} if the functor $\qcoh^{G}(S) \arr \qcoh(S)$ sending $F$ to $F^{G}$ is exact.
\end{definition}


\begin{proposition}
Assume that $S$ is noetherian. Then $G$ is \lr if and only if the functor $\coh^{G}(S) \arr \coh(S)$ defined as $F \mapsto F^{G}$ is exact.
\end{proposition}

\begin{proof}
We need to show that given a surjection $F \arr F'$ of $G$-equivariant \qc sheaves on $S$, the induced morphism $F^{G} \arr {F'}^{G}$ is also surjective. Since $S$ is noetherian,  every \qc sheaf with an action of $G$ is a direct limit of coherent subsheaves with an action of $G$ (see, e.g. Lemma 2.1 \cite{th5} or \cite{Saavedra}). By replacing $F'$ with an arbitrary coherent subsheaf and $F$ with its inverse image in $F$, we may assume that $F'$ is coherent.

Let $\{F_{i}\}$ the inductive system of coherent $G$-equivariant subsheaves of $F$: since $S$ is noetherian and $F'$ is coherent, there will exists some $i$ such that $F_{i}$ surjects onto $F'$. Then $F_{i}^{G} \arr {F'}^{G}$ is surjective, so $F^{G} \arr {F'}^{G}$ is surjective.
\end{proof}

In particular, if $k$ is a field, the category of coherent sheaves on $\spec k$ with an action of $G$ is equivalent to the category of finite-dimensional representations of $G$; hence a finite group scheme over a field is \lr if and only if the functor $V \mapsto V^{G}$, from finite-dimensional representations of $G$ to vector spaces, is exact.

Another, perhaps more customary, way to state this condition is to require that every finite-dimensional representation of $G$ be a sum of irreducible representations.

\begin{proposition}\call{prop:base-change}
Let $S' \arr S$ be a morphism of schemes, $G\arr S$ a group scheme, $G' \eqdef S' \times_{S}G$.

\begin{enumeratea}

\itemref{1} If $G \arr S$ is \lr, then $G' \arr S'$ is \lr.

\itemref{2} If $G' \arr S'$ is \lr and $S' \arr S$ is flat and surjective, then $G \arr S$ is \lr.

\end{enumeratea}
\end{proposition}

\begin{proof}

Let us prove part~\refpart{prop:base-change}{2}. There is a cartesian diagram
   \[
   \xymatrix{
      \cB_{S'}G' \ar[r]^{g}\ar[d]^{\pi'}
   & \cB_{S}G \ar[d]^{\pi}\\
   S' \ar[r]^{f}
   &S
   }
   \]
from which we deduce that the two functors $f^{*}\pi_{*}$ and $\pi'_{*}g^{*}$ are isomorphic. Since $f$ is flat, $g$ is flat as well; also $\pi'_*$ is exact by assumption, so $\pi'_{*}g^{*}$ is exact, hence $f^{*}\pi_{*}$ is exact. But since $f$ is faithfully flat we have that $\pi_{*}$ is exact, as required.

Now for  part~\refpart{prop:base-change}{1}.

First assume that $S'$ is an affine open subscheme of $S$. Then the embedding $j\colon S' \into S$ is quasi-compact, since $S$ is quasi-separated: hence the pushforward $j_{*}$ takes \qc sheaves into \qc sheaves (\cite[I, Corollaire~9.2.2]{EGA}). It is easy to see that if $F'$ is a $G_{S'}$-equivariant \qc sheaf on $S'$, the action of $G_{S'}$ on $F'$ extends to an action of $G$ on $j_{*}F'$: if $T \arr S$ is a flat morphism and $T'$ is the inverse image of $S'$ in $T$, then $G(T')$ acts on $F(T') = j_{*}F(T)$, and this induces an action of $G(T)$ on $j_{*}F(T)$ via the restriction homomorphism $G(T) \arr G(T')$. Then every exact sequence
   \[
   0 \arr F'_{1} \arr F'_{2} \arr F'_{3}\arr 0
   \]
of $G$-equivariant \qc sheaves on $S'$ extends to an exact sequence
   \[
   0 \arr F_{1} \arr F_{2} \arr F_{3}\arr 0
   \]
of $G$-equivariant \qc sheaves on $S$: we take $F_{1} \eqdef j_{*}F'_{1}$ and $F_{2} \eqdef j_{*}F'_{2}$, while $F_{3}$ is defined to be the image of $j_{*}F'_{2}$ in $j_{*}F'_{3}$.
Since taking invariants commutes with restriction to open subschemes, the result follows.

If $\{S_{i}\}$ is an open covering of $S$ by affines, and each restriction $G_{S_{i}}$ is \lr over $S_{i}$, then the disjoint union $\sqcup_{i} G_{S_{i}}$ is \lr over $\sqcup_{i} S_{i}$; we conclude from part~\refpart{prop:base-change}{2} that $G$ is \lr over $S$. Hence being \lr is a local property in the Zariski topology. So to prove part~\refpart{prop:base-change}{1} it suffices to show  that when $S$ and $S'$ are both affine the functor  $ \qcoh^{G'}(S') \arr \qcoh(S')$ is exact.

Then $g$ is also affine, so 
the functor $g_{*}\colon \qcoh^{G'}(S') \arr \qcoh^{G}(S)$ is exact.  By 
assumption $\pi_*$ is exact, therefore  $ \pi_*g_* = f_*\pi'_* $ is exact. But 
the functor $f_*$ has the property that a sequence $F_1 \to F_2 \to F_3$ is 
exact if and only if $f_*F_1 \to f_*F_2 \to f_*F_3$ is exact. It follows that 
$\pi'_*$ is exact, as required.


\end{proof}

\begin{proposition} \call{prop:lr-subgroup}\hfil
The class of \lr group schemes is closed under taking 
\begin{enumeratea}
\itemref{1} subgroup schemes,
\itemref{2} quotients, and 
\itemref{3} extensions. 
\end{enumeratea}
\end{proposition}

\begin{proof}
For part~\refpart{prop:lr-subgroup}{1}, consider a subgroup scheme $G'\subset G$ and the resulting commutative diagram 
$$\xymatrix{
\cB_SG' \ar[r]^i\ar[rd]_{\pi_{G'}} & \cB_SG\ar[d]^{\pi_G} \\& S.
}$$
 It is enough to observe that $i_*= \ind^{G}_{G'}$ is exact (Remark \refall{rmk:important-points}{2}). Since $\pi_{G\,*}$ is exact by assumption, and since $\pi_{G'\,*} \simeq \pi_{G\,*} \circ i_*$ part~\refpart{prop:lr-subgroup}{1} follows. 


For parts~\refpart{prop:lr-subgroup}{2} and~\refpart{prop:lr-subgroup}{3}, consider an exact sequence $$1 \to G' \to G \to G'' \to 1$$ and the corresponding commutative diagram $$\xymatrix{
\cB_SG' \ar[r]^i\ar[rd]_{\pi_{G'}} & \cB_SG\ar[r]^j\ar[d]^{\pi_G} & \cB_SG'' \ar[dl]^{\pi_{G''}}\\& S
}$$

To prove part~\refpart{prop:lr-subgroup}{2}, suppose that $G$ is \lr, so $\pi_{G\,*}$ is exact.  Recall that $j^*$ is exact and $j_*\circ j^*$ is isomorphic to the identity, so $$\pi_{G''\,*} \simeq \pi_{G''\,*}\circ j_*\circ j^*  \simeq \pi_{G\,*}\circ j^*$$ is exact, as required.

For part~\refpart{prop:lr-subgroup}{3}, we have by assumption that $\pi_{G'\,*}$ and $\pi_{G''\,*}$ are exact. Consider the cartesian diagram
$$\xymatrix{
\cB_SG' \ar[r]^i\ar[d]_{\pi_{G'}} &\cB_SG\ar[d]^j  \\ 
  S \ar[r]                            & \cB_SG''.
}$$
The formation of $j_*$ commutes with flat base change on $\cB_SG^{\prime \prime }$ and $S\to \cB_SG''$ is faithfully flat. Thus to verify that $j_*$ is exact it suffices to show that $\pi _{G'*}$ is exact, which holds since $G'$ is assumed linearly reductive.  Therefore $j_*$ is exact (concretely, taking invariants of a $G$ sheaf  by $G'$ is exact even if when consider the induced $G''$-action). So $$\pi_{G\,*} =  \pi_{G''\,*}\circ j_*$$ is exact, as required.
%
%
\end{proof}

\subsection{Classifying \lr group schemes}

We will say that a finite group scheme $\Delta \to S$ is  \emph{diagonalizable} if it is abelian and its Cartier dual is constant. The standard definition only requires the Cartier dual to be constant on the connected components of $S$, (see \cite[Expos\'e VII]{sga3}), but the distinction is of little importance to us. We say that $\Delta \arr S$ is \emph{locally diagonalizable} if its Cartier dual is \'etale.

Given a finite group scheme $\Delta \arr S$, the following conditions are clearly equivalent.

\begin{enumeratea}

\item $\Delta$ is locally diagonalizable.

\item Locally in the \'etale topology, $\Delta$ is diagonalizable.

\item Locally in the fppf topology, $\Delta$ is diagonalizable.

\item Locally in the fpqc topology, $\Delta$ is diagonalizable.

\end{enumeratea}

A finite \'etale group scheme $H \to S$ is said to be \emph{tame} if its degree is prime to all residue characteristics.

\begin{definition}
A group scheme $\pi\colon G \arr S$ is \emph{well-split} if it is isomorphic to a semidirect product $H \ltimes \Delta$, where $H$ is constant and tame and $\Delta$ is diagonalizable.

It is \emph{locally well-split} if there is an fpqc cover $\{S_{i} \arr S\}$, such that the group scheme $S_{i}\times_{S} G  \arr S_{i}$ is well-split for each $i$.
\end{definition}

In characteristic $0$ every finite flat group scheme is \'etale and tame, hence locally constant, hence locally well-split.

\begin{proposition}
Every locally well-split group scheme is \lr.
\end{proposition}

\begin{proof} By  Proposition~\ref{prop:base-change} it suffices to consider well-split group schemes.

The Proposition holds for such group schemes by Proposition \ref{prop:lr-subgroup}~(\ref{prop:lr-subgroup;3}), as it holds for diagonalizable group schemes  (see
\cite[Expos\'e VII]{sga3}), and tame constant group schemes (by Maschke's Lemma).
 \end{proof}

\begin{lemma}\call{lem:over-field} Let $G$ be a locally well-split group scheme over a field $k$, let $\Delta_{0}$ be the connected component of the identity, and $H=G/\Delta_0$.
\begin{enumeratea}

\itemref{1} The group scheme $\Delta_{0}$ is locally diagonalizable, and $H$ is \'etale and tame.

\itemref{2} There exists a finite purely inseparable extension $k'$ of $k$, such that $G_{k'}$ is a semidirect product $H_{k'} \ltimes \Delta_{0\, k'}$. In particular, if $k$ is perfect, then $G$ is a semidirect product $H \ltimes \Delta_{0}$.

\itemref{3} There exists a finite extension $k'$ of $k$ such that $G_{k'}$ is well-split. In particular, if $k$ is algebraically closed, then $G$ is well-split.

\end{enumeratea}

\end{lemma}


\begin{proof}
For part~\refpart{lem:over-field}{1}, we may prove the statement after extending the base field, so we may assume that $G$ is well-split (note that since $\Delta _0$ is a connected scheme with a $k$-rational point, $\Delta _0\otimes _kk'$ is connected for any field extension $k\rightarrow k'$). Then $G$ contains a locally diagonalizable subgroup $\Delta$ such that $G/\Delta$ is \'etale; hence $\Delta_{0}$ coincides with the connected component in $\Delta$; hence it is diagonalizable, and the quotient $H = G/\Delta_{0}$ is \'etale and tame.

For part~\refpart{lem:over-field}{2}, there is a finite purely inseparable extension $k'$ of $k$ such that $(G_{k'})\red$ is smooth over $k'$. Then $(G_{k'})\red$ is a subgroup scheme of $G$, and it maps isomorphically to $H_{k'}$.

For part~\refpart{lem:over-field}{3}, because of part~\refpart{lem:over-field}{2} we may assume that $G$ is of the form $H \ltimes \Delta_{0}$, where $H$ is \'etale and tame and $\Delta_{0}$ is locally diagonalizable. After passing to a finite separable extension of $k$, the group scheme $\Delta_{0}$ becomes diagonalizable and $H$ becomes constant.
\end{proof}

\begin{remark}
In general, if $k$ is not perfect a non-trivial extension is necessary to obtain the splitting. For example suppose $k = k_{0}(a)$, where $a$ is an indeterminate and $k$ of characteristic~$p>2$. Consider the semidirect product $\Gamma \eqdef\rC_{2} \ltimes \mmu_{p}$ over $k$, where $\rC_{2}$ is a cyclic group of order $2$, whose generator $s$ acts on $\mmu_{p}$ as $s \cdot t = t^{-1}$. Conjugation gives a right action of $\mmu_{p}$ on $\Gamma$; given the Kummer extension $k' = k[t]/(t^{p} - a)$, take the quotient $G$ of $\Gamma\times_{k} \spec k'$ by the diagonal action of $\mmu_{p}$. A simple calculation shows that $\mmu_{p}$ acts on the connected component $s\mmu_{p}$ of $\Gamma$  by the formula $y^{-1}(sx)y = sy^{2}x$.

Note that since the conjugation action of $\mmu_p$ on $\Gamma $ is by group homomorphisms, the diagonal action of $\mmu_p$ on $\Gamma \times _{\spec k}\spec k'$ defines descent data for the group scheme $\Gamma \times _{\spec k}\spec k'$ over $\spec k'$, and therefore the quotient $G$ has a group scheme structure.
This group scheme $G$ contains $\mmu_{p}$ as the connected component of the identity, and $G/\mmu_{p}$ is isomorphic to $\rC_{2}$, but the other connected component is a copy of $\spec k[t]/(t^{p} - a^{2})$, whence the sequence is not split.
\end{remark}

\begin{proposition}\label{prop:over-field}
Let $k$ be a field, $G\arr \spec k$ a finite group scheme. Then $G$ is \lr if and only if it is locally well-split.
\end{proposition}

\begin{proof} 
Let $\overline{k}$ be the algebraic closure of $k$; then by Proposition \ref{prop:base-change} the group scheme $G_{\overline{k}}$ is \lr (respectively locally well-split) if and only if $G$ is \lr  (respectively locally well-split); so we may assume that $k$ is algebraically closed. We know that locally well-split groups are \lr, so assume that $G$ is \lr. Denote by $p$ the characteristic of $k$. 

Let $G_{0}$ be the connected component of the identity in $G$. Then $G/G_{0}$ is a \lr constant group. If it were not tame it would contain a subgroup of order $p$, which is not \lr. So we may assume that $G$ is connected, and show that it is diagonalizable.  

The following lemma may be known to the experts, but we have not found a reference.

\begin{lemma}\label{lem:ext-diag}
If a connected finite group scheme $G$ over an algebraically closed field contains a diagonalizable normal subgroup $H$, and $G/H$ is again diagonalizable, then $G$ is also diagonalizable. 
\end{lemma}

\begin{proof}
If we show that $G$ is abelian, then it is diagonalizable: its 
Cartier dual is an extension with a constant quotient and a constant
subgroup,  which is therefore a constant group scheme, because the base field is algebraically closed.

The action by conjugation of $G$ on $H$ defines a homomorphism of
group schemes $G \arr
\underaut_{\mathrm{Gr-Sch}/k}(H)=\underaut_{\mathrm{Gr-Sch}/k}(H^{\rC})$, where $H^{\rC}$ is the
Cartier dual of $H$; but the domain is local, while the target is
constant, so this homomorphism is trivial. Equivalently, $H$ is
central in $G$. 

Let $A$ be a commutative $k$-algebra. The groups $H(A)$ and
$G(A)/H(A)$ are commutative, hence, by ``calculus of commutators"
(see \cite{Gorenstein}, Section 6, in particular Lemma 6.1),  we have a bilinear map  
\begin{align*}
G(A)\times G(A) & \longrightarrow H(A)\\
 (x, y) & \mapsto [x,y]
\end{align*}
This is functorial in $A$, therefore the commutator gives a bilinear
map $G \times G \to H$, and since $H$ is central this gives a bilinear
map $Q \times Q \to H$, where we have set $Q \eqdef G/H$. In particular we get a map of sheaves $Q \to \underhom_{\mathrm{Grp-Sch}/k}(Q,H)$, where both source and target are representable. But again the domain is local and the target is
\'etale, hence  $Q$ is mapped to the trivial map, in other words
the  commutator $Q \times Q \to H$ maps to the identity in $H$. This
means that the commutator is trivial, hence $G$ is abelian.
\end{proof}
%


So we may proceed by induction on the dimension of the vector space
$\H^{0}(G, \cO_{G})$, and assume that $G$ does not contain any nontrivial
normal subgroup schemes. In particular, the Frobenius kernel $G_{1}$ of
$G$ is a normal subgroup scheme of $G$, which does not coincide with
the identity, unless $G$ is trivial: so we have that $G = G_{1}$. In
\cite{Mumford}, p. 139 one says that $G$ has \emph{height
  $1$}. Connected group schemes of height $1$ are classified by their
$p$-Lie algebras (see, e.g., \cite{Mumford}, p. 139). 

\begin{lemma}[Jacobson \cite{Jacobson}, Chapter 5, Exercise 14, p. 196]
Let $G$ be a non-abelian group scheme of height 1. Then $G$ contains $\aa_p$, and hence is not \lr.
\end{lemma}

\begin{proof}
Considering the $p$-lie algebra $\frg$ of $G$, we need to find an
element $w\in \frg$ such that $w^{p} = 0$. Since $\frg$ is finite
dimensional, for each $v\in \frg$ there is a minimal $n$ such that
$\{v, v^{p}, v^{p^2},\ldots,v^{p^n}\}$ is linearly dependent,
giving a monic $p$-polynomial $$f_v(x) = x^{p^n} + a^{(v)}_{n-1}x^{p^{n-1}} + \cdots +
a^{(v)}_0 x$$  such that $f_v(v) = 0$. 

Note that if $a^{(v)}_0 = 0$  
then the nonzero element $$w  \ \eqdef\  f_v^{1/p} (v)\ \  =\ \  v^{p^{n-1}}\ +
\ (a^{(v)}_{n-1})^{1/p}\, v^{p^{n-2}}\ + \cdots + \
(a^{(v)}_1)^{1/p}\, v$$ satisfies $w^p = 0$. So, arguing by
contradiction,  we may assume that
$a^{(v)}_0 \neq 0$ for 
all nonzero $v$, i.e. $f_v$ is separable. Since the minimal polynomial of
$\operatorname{ad}(v)$ divides  $f_v$, we have that
$\operatorname{ad}(v)$ is semisimple for every nonzero $v$. 

Since
$\frg$ is assumed non-commutative, there is $v$ with
$\operatorname{ad}(v)\neq 0$, hence it has a nonzero eigenvector $v'$
with nonzero eigenvalue. But then the action of
$\operatorname{ad}(v')$ on $\operatorname{Span}(v,v')$ is nonzero and 
nilpotent, contradicting semisimplicity.
\end{proof} 

Back to the proposition, we deduce that $G$ is abelian. Every subgroup
scheme is normal, so $G$ cannot contain any proper subgroup scheme. But by Cartier duality the only local abelian group schemes with this property are $\aa_{p}$ and $\mmu_{p}$; and again $\aa_{p}$ is not \lr. Hence $G= \mmu_{p}$, and we are done. 
\end{proof}

This completes our analysis for the case of group schemes over fields. To handle the general case we need the following fact.

\begin{lemma}\label{lem:fiber->local}
Let $G \arr S$ be a finite flat group scheme of finite presentation. Assume that there is a point $s = \spec k(s) \in S$ such that the
fiber $G_{s} \arr \spec k(s)$ is locally well-split. Then there exists
a flat quasi-finite map $U \arr S$ of finite presentation, whose image
includes $s$, such that $G_{U}$ is well-split. 

In particular, let $V$ be the image of $U$ in $S$, which is open; then
the restriction $G_{V} \arr V$ is locally well-split. 
\end{lemma}

\begin{proof}
By standard arguments, we may assume that $S$ is connected, affine and
of finite type over $\ZZ$. By Lemma~\refall{lem:over-field}{3}, there is a finite extension $k$ of $k(s)$
such that $G_{k}$ is of the form $H \ltimes\Delta_{0}$, where
$\Delta_{0}$ is a connected diagonalizable group scheme and $H$ is a constant group scheme, associated with a finite group $\Gamma$. After
base change by a quasi-finite flat morphism over $S$, we may assume that
$k(s) = k$. The group scheme $\Delta_{0}$ extends uniquely to a
diagonalizable group scheme $\Delta_{0}$ on $S$, that we still denote
by $\Delta_{0}$. Also, we denote again by $H$ the group scheme over
$S$ associated with $\Gamma$; the action of $H$ on $\Delta_{0}$, which is
defined over $s$, extends uniquely to an action of $H$ on
$\Delta_{0}$. Set $G' = H \ltimes \Delta_{0}$. We claim that $G$ and
$G'$ become isomorphic after passing to a flat morphism of finite type
$U \arr S$, whose image includes $s$. We use the following lemma, which shows that after passing to an \'etale neighborhood of $s$ in $S$ 	there exists a $(G, G')$-bitorsor $I \arr S$. Given such a bitorsor, we have that $G$ is the group scheme of automorphisms of $I$ as a $G'$-torsor, and is thus the twisted form of $G'$ coming from $I$ and from the homomorphism $G' \arr \aut(G')$ given by conjugation; so if $I \arr S$ has a section, $G$ and $G'$ are isomorphic. Hence the pullbacks of $G$ and $G'$ to $I$ are isomorphic, 
 and $I \arr S$ is flat and quasi-finite. This completes the proof.
\end{proof}

\begin{lemma}
Let $G \arr S$ and $G' \arr S$ be two group schemes over $S$, and let $s \in S$ be a point such that the fibers $G_{s}$ and $G'_{s}$ are \lr and isomorphic. Then there exists an \'etale neighborhood $s \arr U \arr S$ of $s$ in $S$ and a $(G_{U}, G'_{U})$-bitorsor $I \arr U$.
\end{lemma}

\begin{proof}
We can pass to the henselization $R$ of the local ring $\cO_{S,s}$, and assume that $S = \spec R$ and that $s$ is the closed point of $S$. Call $\frm$ the maximal ideal of $R$ and set $R_{n} \eqdef R/\frm^{n+1}$ and $S_{n} \eqdef \spec R_{n}$. We will show that there exist a sequence of $(G,G')$-bitorsors $I_{n} \arr S_{n}$, such that the restriction of each $I_{n}$ to $S_{n-1}$ is isomorphic to $I_{n-1}$; then the result follows from Artin's approximation theorem.

We present two methods of proof, one abstract and one more explicit. Both use deformation theory.

{\sc Method 1: rigidity using the cotangent complex.}

By \cite[Remarque~1.6.7]{giraud}, a $(G,G')$-bitorsor $I_{n} \arr S_{n}$ is the same as an equivalence of fibered categories $\rho _n:\cB_{S_{n}}G_{S_{n}} \simeq \cB_{S_{n}}G'_{S_{n}}$. We inductively construct a compatible system of such isomorphisms.  So suppose $\rho _n$ has been constructed.  We then wish to find a dotted arrow filling in the diagram
\begin{equation}\label{liftrho}
\xymatrix{\cB_{S_n}G_n\ar[d]^{\rho _n}\ar@{^{(}->}[r]& \cB _{S_{n+1}}G_{n+1}\ar@{.>}[d]^{\rho _{n+1}}\\
\cB_{S_n}G_n'\ar@{^{(}->}[r]& \cB_{S_{n+1}}G_{n+1}'.}
\end{equation}
Note that any such morphism $\rho _{n+1}$ is automatically an isomorphism since $G_{n+1}$ and $G_{n+1}'$ are flat over $S_{n+1}$.   Let $\bL_{\cB G_s}\in D_{\text{coh}}(\mc O_{\cB G_s})$ denote the cotangent complex of $\cB G_s$ over $s$ as defined in \cite{L-MB} (and corrected in \cite{Olsson-sheaves}).

\begin{lem} We have $\bL_{\cB G_s}\in D^{[0,1]}_{\text{\rm coh}}(\mc O_{\cB G_s})$
\end{lem}
\begin{proof} 
Consider the map $p:s\rightarrow \cB G_s$ corresponding to the trivial torsor.  The map $p$ is faithfully flat so it suffices to show that $p^*\bL_{\cB G_s}$ has cohomology concentrated in degrees $0$ and $1$.  From the distinguished triangle
$$
p^*\bL _{\cB G_s}\rightarrow \bL_{s/s}\rightarrow \bL_{s/\cB G_s}\rightarrow p^*\bL _{\cB G_s}[1]
$$
and the fact that $\bL _{s/s} = 0$, we see that $p^*\bL _{\cB G_s}\simeq \bL _{s/\cB G_s}[-1]$.  Therefore it suffices to show that $\bL _{s/\cB G_s}$ is concentrated in degrees $-1$ and $0$.  Consider the cartesian diagram
$$
\xymatrix{s\ar[d]_{p}& G_s\ar[l]_q \ar[d]\\
\cB G_s& s\ar[l]_-p.}
$$
Since $q$ is faithfully flat it then suffices to show that
$$
q^*\bL _{s/\cB G_s}\simeq \bL _{G_s/s}
$$
is concentrated in degrees $-1$ and $0$.   This follows from the fact that $G_s$ is a local complete intersection (this can be seen for example from \ref{prop:over-field}), and \cite[III.3.2.6]{Illusie}.
\end{proof}

Since $s$ is the spectrum of a field any coherent sheaf on $\cB G_s$ is locally free, and therefore for any coherent sheaf $\mc F$ on $\cB G_s$ we have
$$
\mls R\mls Hom(\bL _{\cB G_s}, \mc F)\in D_{\text{coh}}^{[-1, 0]}(\mc O_{\cB G_s}).
$$
Since the global section functor is exact on the category $Coh(\mc O_{\cB G_s})$ (since $G_s$ is linearly reductive) we obtain that
$$
\text{Ext}^i(\bL _{\cB G_s}, \mc F) = 0, \ \ \text{for $i\neq -1,0$}.
$$

By \cite[1.5]{deform}, the obstruction to finding the arrow $\rho _{n+1}$ filling in \ref{liftrho}, is a class in the group
$$
\text{Ext}^1(\bL _{\cB G_s}, \mathfrak{m}^{n+1}/\mathfrak{m}^{n+2}\otimes _{k(s)}\mc O_{\cB G_s})
$$
which by the above is zero.  It follows that there exists an arrow $\rho _{n+1}$ filling in \ref{liftrho}.

{\sc Method 2: lifting using Lie algebras.}

Set $G_{n} \eqdef G_{S_{n}}$ and $G'_{n} \eqdef
G'_{S_{n}}$. 

Let us start from the tautological $G_{0}$-torsor $s = S_{0} \arr
\cB_{S_{0}}G_{0}$, which we think of as a $G'$-torsor via an isomorphism $G_{0}  \simeq G'_{0}$. Our aim now is
to construct a sequence of $G'$-torsors $P_{n} \arr \cB_{S_{n}}G_{n}$,
such that the restriction of each $P_{n}$ to $S_{n-1}$ is isomorphic
to the $G'$-torsor $P_{n-1} \arr \cB_{S_{n-1}}G_{n-1}$. 

The Lie algebra $\frg$ of $G_{0} = G'_{0}$ is a representation of
$G_{0}$, corresponding to a coherent sheaf on $\cB_{S_{0}}G_{0}$. It
is well known that the obstruction to extending $P_{n-1} \arr
\cB_{S_{n-1}}G_{n-1}$ to a $G'$-torsor lies in the sheaf cohomology
$\H^{2}\bigl(\cB_{S_{0}}G_{0}, (\frm^{n}/\frm^{n+1}) \otimes
\frg\bigr)$; and this coincides with the cohomology of $G_{0}$ in the
representation $(\frm^{n}/\frm^{n+1}) \otimes \frg$, which is $0$,
because $G_{0}$ is \lr. 

Each $G'$-torsor $P_{n} \arr \cB_{S_{n}}G_{n}$ yields a $(G,G')$-bitorsor
   \[
   I_{n} \eqdef S_{n}\times_{\cB_{S_{n}}G_{n}} P_{n} \arr S_{n},
   \]
where the morphism $S_{n} \arr \cB_{S_{n}}G_{n}$ is the one given by
the trivial torsor $G_{n} \arr S_{n}$. So we obtain the desired sequence of bitorsors.
\end{proof}

Here is our main result on \lr group schemes.

\begin{theorem}\label{Th:main-lr} Let $G \arr S$ be a finite flat group scheme of finite presentation. The following conditions are equivalent. 

\begin{enumeratea}

\item $G \arr S$ is \lr.

\item $G \arr S$ is locally well-split.

\item The fibers of $G \arr S$ are \lr.

\item The geometric fibers of $G \arr S$ are well-split.

\end{enumeratea}

Furthermore, if $S$ is noetherian these conditions are equivalent to either of the following two conditions.

\begin{enumeratea}\setcounter{enumi}{4}

\item The closed fibers of $G \arr S$ are \lr.

\item The geometric closed fibers of $G \arr S$ are well-split.

\end{enumeratea}

\end{theorem}

\begin{proof}
This follows from Proposition~\ref{prop:over-field}, Lemma~\ref{lem:over-field} and Lemma~\ref{lem:fiber->local}.
\end{proof}

For later use we prove the following.

\begin{lemma}\label{lem:sh->wellsplit}
Let $G \arr \spec R$ be a \lr group scheme, where $R$ is a henselian local ring. Then there exists a locally diagonalizable normal subgroup $\Delta$ of $G$ such that $G/\Delta$ is \'etale and tame.

If $R$ is strictly henselian, then $\Delta$ is diagonalizable and $G/\Delta$
is constant.
\end{lemma}

\begin{proof}
We may assume that $R$ is the  henselization of a scheme of finite type over $\ZZ$ at a point. Set $S = \spec R$, and let $\frm$ and $k$ be the maximal ideal and the residue field of $R$. Let $\Delta_{0}$ be the connected component of the identity in $G_{k}$; by Lemma~\refall{lem:over-field}{1}, it is locally diagonalizable. Then we claim that there exists a unique closed subscheme $\Delta \subseteq G$, flat over $S$, whose restriction to $G_{k}$ coincides with $\Delta_{0}$. If $R$ is artinian then $G$ and $G_{k}$ are homeomorphic, and $\Delta$ is simply the connected component of the identity in $G$, which is automatically flat. In general, call $\Delta_{n}$ the connected component of the identity in $G_{S_{n}}$, where $S_{n} \eqdef \spec R/\frm^{n+1}$. If $\Delta$ is such a subscheme of $G$, we have $\Delta \cap \Delta_{n} = \Delta_{n}$ for each $n$; hence $\Delta$ is unique. Artin's approximation theorem, applied to the functor that sends each $R$-algebra $A$ into the set of subschemes of $G_{A}$ that are flat over $A$, ensures the existence of such a $\Delta$.

Now we need to show that it is a subgroup scheme, which is equivalent to showing that it is closed under multiplication and inverses. Let $m\colon G \times G \arr G$ be the multiplication map; $m^{-1}(\Delta) \cap (G_{S_{n}} \times G_{S_{n}})$ contains $(\Delta \times \Delta) \cap (G\times G)_{S_{n}} = \Delta_{n} \times \Delta_{n}$ for each $n$, because $\Delta_{n}$ is a subgroup scheme of $G_{S_{n}}$; hence $\Delta \times \Delta \subseteq m^{-1}(\Delta)$ as required. The argument for inverses is similar.

The group scheme of automorphisms of $\Delta_n$ over $S_n$ is unramified, since its closed fiber is \'etale. It follows that the image of $\Delta_n$ in its automorphism group scheme is trivial, hence $\Delta_n$ is abelian.
 
It is easy to see that each $\Delta_{n}$ is locally diagonalizable, by looking at Cartier duals, and that $G_{n}/\Delta_{n}$ is tame and \'etale. From this it follows that $\Delta$ is locally diagonalizable and $G/\Delta$ is tame and \'etale.

The last statement follows from the fact that any \'etale cover of the spectrum of a strictly henselian ring is trivial.
\end{proof}

\subsection{A remark on group cohomology}\label{Sec:group-cohomology}
The following is a completely elementary fact, which we explain for lack of a suitable reference.

Let
   \[
   1 \arr A \arr E \arr G \arr 1
   \]
be an extension of groups, with $A$ abelian. Consider the induced action of $G$ on $A$ by conjugation. Conjugation by an element of $A$ gives an automorphism of $E$, which induces the identity on both $A$ and $G$.

If $\phi\colon E \arr E$ is another such automorphism, we can consider the function $E \arr A$ defined by $u \mapsto \phi(u)u^{-1}$; this descends to a function $\psi\colon G \arr A$, linked with $\phi$ by the formula $\phi(u) = \psi([u])u$, where $[u]\in G$ denotes the image of $u$. It is easy to see that $\psi$ is a crossed homomorphism, and that sending $\phi$ into $\psi$ gives an isomorphism of the group of automorphisms of $E$ that induce the identity on $A$ and on $G$ with the group $\rZ^{1}(G, A)$ of crossed homomorphisms. It is also easy to see that $\phi$ is given by conjugation by an element of $A$ if and only if the crossed homomorphism is a boundary. Hence if $\H^{1}(G,A) = 0$ the only such automorphisms are obtained by conjugating by elements of $A$.

\subsection{\'Etale local liftings of \lr group schemes}

We will need the following result.

\begin{proposition}\label{prop:extend-lr}
Let $S$ be a scheme, $p \in S$ a point, $G_{0} \arr p$ a \lr group scheme. There exists an \'etale morphism $U \arr S$, with a point $q \in U$ mapping to $p$, and a \lr group scheme $\Gamma \arr U$ whose restriction $\Gamma_{q} \arr q$ is isomorphic to the pullback of $G_{0}$ to $q$.
\end{proposition}

Let us start with a lemma. Let $k$ be a field, $G \arr \spec k$ a well-split group scheme. Let $\Delta$ be the connected component of the identity of $G$, $H = G/\Delta$. Call $\underaut_{k}(G)$ the group scheme representing the functor of automorphisms of $G$ as a group scheme: there is a homomorphism $\Delta \arr \underaut_{k}(G)$ sending each section of $\Delta$ into the corresponding inner automorphism of $G$; this induces an embedding $\Delta/\Delta^{H} \subseteq \underaut_{k}(G)$, where $\Delta ^H$ denotes the $H$-invariants of $\Delta $.

\begin{lemma}\label{lem:inner-autos}
The connected component of the identity of $\underaut_{k}(G)$ is $\Delta/\Delta^{H}$.
\end{lemma}

\begin{proof}
Since $\Delta$ is a characteristic subgroup scheme of $G$, each automorphism of $G_{A} \arr \spec A$, where $A$ is a $k$-algebra, preserves $\Delta_{A}$. Hence we get homomorphisms of group schemes $\underaut_{k}(G) \arr \underaut_{k}(\Delta)$ and $\underaut_{k}(G) \arr \underaut_{k}\bigl(H\bigr)$, inducing a homomorphism
   \[
   \underaut_{k}(G) \arr \underaut_{k}(\Delta)\times \underaut_{k}\bigl(H\bigr);
   \]
the kernel of this homomorphism contains $\Delta/\Delta^{H}$. Let us denote by $E$ this kernel; since $\underaut_{k}(\Delta)\times \underaut_{k}\bigl(H\bigr)$ is \'etale over $\spec k$, it is enough to prove that $E$ coincides with $\Delta/\Delta^{H}$.

To do this, we may pass to the algebraic closure of $k$, and assume that $k$ is algebraically closed; then it is enough to prove that given a $k$-algebra $A$, for any element $\alpha \in E(A)$ there exists a faithfully flat extension $A \subseteq A'$ such that the image of $\alpha$ in $E(A')$ comes from $(\Delta/\Delta^{H})(A')$.

By passing to a faithfully flat extension, we may assume that $G(B) \arr H(B)$ is surjective for any $A$-algebra $B$ (because $H$ is constant), so we have an exact sequence
   \[
   1 \arr \Delta(B) \arr G(B) \arr H(B) \arr 1.
   \]
Furthermore, again because $H$ is constant, for any $A$-algebra $B$ we have
   \[
   \Delta^{H}(B) = \Delta(B)^{H(B)};
   \]
hence for any $B$ we have an injective homomorphism
   \[
   \Delta(B)/\Delta(B)^{H(B)} \arr (\Delta/\Delta^{H})(B).
   \]
Let us show that $\alpha$ comes from $(\Delta/\Delta^{H})(A)$.

Set $B = \Gamma(G_{A}, \cO)$, so that $G_{A} = \spec B \arr \spec A$. Then it is easy to see that the natural restriction homomorphism $\aut_{A}(G_{A}) \arr \aut\bigl(G(B)\bigr)$ is injective. The group $\Delta(B)$ has an order that is a power of the characteristic of $k$, while the order of $H(B)$ is prime to the characteristic; so $\H^{1}\bigl(H(B), \Delta(B)\bigr) = 0$. By the discussion in \ref{Sec:group-cohomology}  there exists an element $\delta_{B}$ of $\Delta(B)$ whose image in $\aut\bigl(G(B)\bigr)$ (the automorphism given by conjugation by $\delta _B$) coincides with the image of $\alpha$.

Call $\overline{\delta}_{B}$ the image of $\delta_{B}$ in $(\Delta/\Delta^{H})(B)$. We claim that $\overline{\delta}_{B}$ is the image of an element $\overline{\delta}$ of $(\Delta/\Delta^{H})(A)$; then the image of $\delta$ in $E(A)$ must be $\alpha$, because $\aut_{A}(G_{A})$ injects into $\aut\bigl(G(B)\bigr)$.

To prove this, since $(\Delta/\Delta^{H})(A)$ is the equalizer of the two natural maps $(\Delta/\Delta^{H})(B) \double (\Delta/\Delta^{H})(B \otimes_{A} B)$, it is enough to show that the two images of $\overline{\delta}_{B}$ in $(\Delta/\Delta^{H})(B \otimes_{A} B)$ coincide. The two images of $\delta_{B}$ in $\aut\bigl(G(B \times_{A} B)\bigr)$ are equal; since $\Delta(B \times_{A} B)/\Delta(B \times_{A} B)^{H(B \times_{A} B)}$ injects into $\aut\bigl(G(B \times_{A} B)\bigr)$, this implies that the two images of $\delta_{B}$ into $\Delta(B \times_{A} B)^{H(B \times_{A} B)}$ coincide. The images of these via the natural injective homomorphism
   \[
   \Delta(B \times_{A} B)/\Delta(B \times_{A} B)^{H(B \times_{A} B)} \arr
   (\Delta/\Delta^{H})(B \times_{A} B)
   \]
are the two images of $\overline{\delta}_{B}$, and this completes the proof.
\end{proof}

\begin{proof}[Proof of Proposition~\ref{prop:extend-lr}] Let $\overline{k(p)}$ be the algebraic closure of $k(p)$; the pullback $G_{\overline{k(p)}}$ is well-split, that is, it is a semi-direct product $H_{\overline{k(p)}} \ltimes \Delta_{\overline{k(p)}}$, where $H_{\overline{k(p)}}$ is \'etale, hence a constant group, and $\Delta_{\overline{k(p)}}$ is connected and diagonalizable. This is the pullback of a group scheme $\Gamma = H \ltimes \Delta \arr S$, where $H$ is constant and $\Delta$ is diagonalizable; passing to a Zariski open neighborhood of the image of $p$ in $S$, we may assume that $H$ is tame, so $\Gamma$ is well split. The group scheme $G_{0}$ is a twisted form of the fiber $\Gamma_{p}$. So we need to show that every twisted form of $\Gamma_{p} \arr S$ on a point $p \in S$ extends to an \'etale neighborhood of $p$. This twisted form is classified by an element of the non-abelian cohomology group $\H^{1}\fppf\bigl(p, \underaut_{k}(\Gamma)\bigr)$. Let us set $\Delta' = \Delta/\Delta^{H}$. The quotient $\underaut_{k}(\Gamma)/\Delta'_{p}$ is \'etale, by Lemma~\ref{lem:inner-autos}; hence the image of this element in
   \[
   \H^{1}\fppf\bigl(p, \underaut_{k}(\Gamma)/\Delta'_{p}\bigr) = 
   \H^{1}\et\bigl(p, \underaut_{k}(\Gamma)/\Delta'_{p}\bigr)
   \]
is killed after passing to a finite separable extension of $k(p)$. Any such extension is of the form $k(q)$, where $U \arr S$ is an \'etale map and $q$ is a point on $U$ mapping on $p$. We can substitute $S$ with $U$, and assume that the image of our element of $\H\fppf\bigl(p, \underaut_{k}(\Gamma)\bigr)$ in $\H^{1}\fppf\bigl(p, \underaut_{k}(\Gamma)/\Delta'_{p})$ is trivial. We have an exact sequence of pointed sets
   \begin{align*}
   \H^{1}\fppf\bigl(p, \Delta'_{p}\bigr) \arr
   \H^{1}\fppf\bigl(p, \underaut_{k}(\Gamma)\bigr) \arr
   \H^{1}\fppf\bigl(p, \underaut_{k}(\Gamma)/\Delta'_{p}\bigr);
   \end{align*}
so our element comes from $\H^{1}\fppf(p, \Delta'_{p})$. Since $\Delta'$ is diagonalizable, it is enough to prove that every element of $\H^{1}\fppf(p, \mmu_{n})$ comes from $\H^{1}\fppf(S, \mmu_{n})$, after restricting $S$ in the Zariski topology.

By Kummer theory, every $\mmu_{n}$-torsor over $k(p)$ is of the form
   \[
   \spec k(p)[t]/(t^{n} - a) \arr \spec k(p)
   \]
for some $a \in k(p)^{*}$, with the obvious action of $\mmu_{n}$ on  $\spec k(p)[t]/(t^{n} - a)$. After passing to a Zariski neighborhood of $p \in S$, we may assume that $a$ is the restriction of a section $f \in \cO^{*}(S)$. Then the $\mmu_{n}$-torsor $\spec_{S} \cO_{S}[t]/(t^{n} - f) \arr S$ restricts to $\spec k(p)[t]/(t^{n} - a) \arr \spec k(p)$, and this completes the proof.
\end{proof}

\section{Tame stacks}\label{section3}

Let $S$ be a scheme, $\cM \arr S$ a locally finitely presented algebraic stack
over $S$. We denote by $\cI \arr \cM$ the inertia group stack; we will
always assume that $\cI \arr \cM$ is finite (and we say that $\cM$ has
\emph{finite inertia}). If $T \arr S$ is a morphism, and $\xi$ is an
object of $\cM(T)$, then the group scheme $\underaut_{T}(\xi) \arr T$
is the pullback of $\cI$ along the morphism $T \arr \cM$ corresponding
to $\xi$. 

Under this hypothesis, it follows from \cite{Keel-Mori} that there exists a moduli space $\rho\colon \cM \arr M$; the morphism $\rho$ is proper.

\begin{definition}\label{Def:tame}
The stack $\cM$ is \emph{tame} if the functor $\rho_{*}\colon \qcoh{\cM} \arr \qcoh{M}$ is exact.
\end{definition}

When $G \arr S$ is a finite flat group scheme, then the moduli space of $\cB_{S}G \arr S$ is $S$ itself; so $\cB_{S}G$ is tame if and only if $G$ is \lr.

\begin{theorem}\call{thm:main}
The  following conditions are equivalent.
\begin{enumeratea}

\itemref{1} $\cM$ is tame.

\itemref{2} If $k$ is an algebraically closed field with a morphism $\spec k \arr S$ and $\xi$ is an object of $\cM(\spec k)$, then the automorphism group scheme $\underaut_{k}(\xi) \arr \spec k$ is \lr.

\itemref{3} There exists an fppf cover $M'\arr M$, a \lr group scheme $G \arr M'$ acting an a finite and finitely presented scheme $U \arr M'$, together with an isomorphism
   \[
   \cM \times_{M} M' \simeq [U/G] 
   \]
of algebraic stacks over $M'$.

\itemref{4} Same as \refpart{thm:main}{3}, but $M' \arr M$ is assumed to be \'etale and surjective.

\end{enumeratea}
\end{theorem}

\begin{corollary}\call{2.3} Let $\mc M$ be a tame stack over a scheme $S$ and let $\mc M\rightarrow M$ be its moduli space.

\begin{enumeratea}

\itemref{1} If $M'\arr M$ is a morphism of algebraic spaces, then the moduli space of $M'\times _M\cM$ is $M'$.

\itemref{2} If $\cM$ is flat over $S$ then $M$ is also flat over $S$.

\end{enumeratea}
\end{corollary}

\begin{proof}
This corollary is proved by standard arguments as follows.

Formation of moduli spaces commutes with flat base change; furthermore, if $\cM$ is an algebraic stack locally of finite presentation with finite diagonal, $\cM \arr M$ is a morphism into an algebraic space, and $\{M_{i} \to M\}$ is an fppf cover, then $M$ is the moduli space of $\cM$ if and only if for each $i$ the algebraic space $M_{i}$ is the moduli space of $M_{i} \times_{M} \cM$.

Hence by Theorem~\ref{thm:main} we may assume that $S = \spec R$ and $M = \spec A$ are affine, and that $\cM$ is of the form $[U/G]$, where $G$ is a \lr group scheme. Then $U$ is finite over $M = U/G$, so it is affine: write $U = \spec B$. Then the group scheme $G$ acts on the $A$-algebra $B$, and $A = B^{G}$.

Let $N$ be an $R$-module; then there is a natural action of $G$ on $N\otimes_{R} B$, obtained from the action of $G$ on $B$. There is a natural homomorphism $N \otimes_{R} A \arr N\otimes_{R} B$, which is $G$-equivariant, when letting $G$ act trivially on $N \otimes_{R} A$. We claim that the induced homomorphism
   \[
   N \otimes_{R} A \arr (N \otimes_{R} B)^{G}
   \]
is an isomorphism. This is obvious when $N$ is free, because in this case $N\otimes_{R} B$ is a direct sum of copies of $B$, with $G$ acting separately on each copy. In general, let $F_{1} \arr F_{0} \arr N \arr 0$ be a free presentation of $N$. We have a commutative diagram
   \[
   \xymatrix{
   F_{1}\otimes_{R}A \ar[r]\ar[d] & F_{0}\otimes_{R}A \ar[r]\ar[d]
      &N \otimes_{R} A \ar[r]\ar[d] & 0\\
   (F_{1}\otimes_{R}B)^{G} \ar[r] &(F_{0}\otimes_{R}B)^{G} \ar[r]
      &(N \otimes_{R} B)^{G} \ar[r] & 0\\      
   }
   \]
in which both rows are exact (the second one because $G$ is \lr). The first two columns are isomorphisms, hence so is the third.

Part~\refpart{2.3}{1} follows easily from this. We need to check that if $M' \arr M$ is a morphism, the projection $M' \times_{M} \cM  \arr M'$ makes $M'$ into the moduli space or $M' \times_{M} \cM$. We may assume that $M' = \spec A'$ is affine, and that $S = M$, so that $R = M$; then the statement follows from the fact the natural homomorphism $A' = A' \otimes_{A} A \arr (A' \otimes_{A} B)^{G}$ is an isomorphism. 

Part~\refpart{2.3}{2} is also easy. Assume that $B$ is flat over $R$. The isomorphism $N \otimes_{R} A \simeq (N \otimes_{R} B)^{G}$ is functorial in the $R$-module $N$, and both tensoring with $B$ and taking invariants give exact functors.
\end{proof}

\begin{corollary}
If $\cM \arr S$ tame and $S' \arr S$ is a morphism of schemes, then $S' \times_{S} \cM$ is a tame stack over $S'$.
\end{corollary}

\begin{corollary}
The stack $\cM \arr S$ tame if and only if for any morphism $\spec k \arr S$, where $k$ is an algebraically closed field, the geometric fiber $\spec k\times_{S}\cM$ is tame.
\end{corollary}

\begin{proof}[Proof of Theorem~\ref{thm:main}] It is obvious that \refpart{thm:main}{4} implies \refpart{thm:main}{3}. It is straightforward to see that \refpart{thm:main}{3} implies both \refpart{thm:main}{1} and \refpart{thm:main}{2}.

Let us check that \refpart{thm:main}{1} implies \refpart{thm:main}{2}. Let $\spec k  \arr \cM$ be the morphism corresponding to the object $\xi$ of $\cM(\spec k)$; set $G = \underaut_{k}(\xi)$. Call $\cM_{0}$ the pullback $\spec k \times_{M} \cM$; this admits a section $\spec k \arr \cM_{0}$, and the residual gerbe of this section, which is a closed substack of $\cM_{0}$, is isomorphic to $\cB_{k}G$. So we get a commutative (non cartesian) diagram
   \[
   \xymatrix{
   \cB_{k}G \ar[d]^{\rho'}\ar[r]^-{g}
   & \cM\ar[d]^{\rho} \\
   \spec k \ar[r]^-{f} 
   & M
   }
   \]
whose rows are affine. Hence we have that $g_{*}\colon \qcoh(\cB_{k}G) \arr \qcoh(\cM)$ is an exact functor, while $\rho_{*}\colon \qcoh(\cM) \arr \qcoh(M)$ is  exact by definition. Also we have an equality of functors $f_{*}\rho'_{*} = \rho_{*}g_{*}$; hence, if
   \[
   0 \arr V_{1} \arr V_{2} \arr V_{3} \arr 0
   \]
is an exact sequence of representations of $G$, considered as an exact sequence of \qc sheaves on $\cB_{k}G$, we have that the sequence
   \[
   0 \arr f_{*}\bigl(V_{1}^{G}\bigr) \arr f_{*}\bigl(V_{2}^{G}\bigr)
   \arr f_{*}\bigl(V_{3}^{G}\bigr) \arr 0
   \]
is exact; and this implies that
   \[
   0 \arr V_{1}^{G} \arr V_{2}^{G} \arr V_{3}^{G} \arr 0
   \]
is exact. Hence $G$ is \lr, as claimed.

Now let us prove that \refpart{thm:main}{2} implies \refpart{thm:main}{4}. In fact, we will prove a stronger version of this implication.

\begin{proposition}\label{prop:stronger-version}
Let $\cM \arr S$ be a locally finitely presented algebraic stack with finite inertia and moduli space $\rho\colon \cM \arr M$. Let $k$ be a field with a morphism $\spec k \arr S$, and let $\xi$ be an object of $\cM(\spec k)$; assume that the automorphism group scheme $\underaut_{k}(\xi)\arr \spec k$ is \lr. Denote by $p\in M$ the image of the composite $\spec k \arr \cM \arr M$. Then there exists an \'etale morphism $U \arr M$ having $p$ in its image, a \lr group scheme $G \arr U$ acting on a finite scheme $V \arr U$ of finite presentation, and an isomorphism $[V/G] \simeq U \times_{S}\cM$ of algebraic stacks over $U$.
\end{proposition}

Thus, if $\cM$ has an object over a field with \lr automorphism group, then there is an open tame substack of $\cM$ (the image of $U$) containing this object.

\begin{proof}
The proof is divided into three steps.

We may assume that $M$ is affine and of finite type over $\ZZ$.

{\sc The case $k = k(p)$.\enspace} We start by assuming that the residue field $k(p)$ of $p \in M$ equals $k$. After passing to an \'etale cover of $M$, we may also assume that $\underaut_{k}(\xi)$ extends to a \lr group scheme $G \arr M$ (Proposition~\ref{prop:extend-lr}).

By standard limit arguments we may assume that $M$ is the spectrum of a local henselian ring $R$ with residue field $k$. The result follows once we have shown that there is  a \emph{representable} morphism $\cM
\arr \cB_{M}G$ of algebraic stacks, equivalently, a $G$-torsor $P
\arr \cM$ in which the total space is an algebraic space.

Let us denote by $\cM_{0}$ the residual gerbe
$\cB_{k}\underaut_{k}(\xi) = \cB_{k}G_{p}$; this is a closed substack
of $\cM$, having $\spec k$ as its moduli space. This closed substack
gives a sheaf of ideals $\cI \subseteq \cO_{\cM}$; we denote by
$\cM_{n}$ the closed substack of $\cM$ whose sheaf of ideals is
$\cI^{n+1}$. Denote by $\frg$ the Lie algebra of $\underaut_{k}(\xi)$.

The obstruction to extending a $G$-torsor $P_{n-1} \arr \cM_{n-1}$ to a $G$-torsor $P_{n} \arr \cM_{n}$ lies in
   \begin{align*}
   \H^{2}\bigl(\cM_{0}, (\cI^{n}/\cI^{n+1})\otimes \frg\bigr) &=
   \H^{2}\bigl(G_{p}, (\cI^{n}/\cI^{n+1})\otimes \frg\bigr)\\
   & = 0;
   \end{align*}
Alternatively, in terms of the cotangent complex, the obstruction lies
in $\ext^1(\bL g^*\bL_{\cB_kG_p},\cI^{n}/\cI^{n+1}) =0$, where $g:\mc M_0\rightarrow \cB _MG$ is the morphism defined by $P_0$.

Hence we can construct a sequence of $G$-torsors $P_{n} \arr \cM_{n}$,
such that the restriction of $P_{n}$ to $P_{n-1}$ is isomorphic to
$P_{n-1}$, and such that the torsor $P_{0} \arr\cM_{0}$ has $\spec k$
as its total space. 

Let $\frm$ be the maximal ideal of $R$, and
set $M_{n} = \spec R/\frm^{n+1}$. The systems of ideals $\{\cI^{n}\}$ and
$\frm^{n}\cO_{\cM}$ are cofinal; hence we get a sequence of
$G$-torsors $Q_{n} \arr M_{n}\times_{M}\cM$, such that the restriction
of $Q_{n}$ to $Q_{n-1}$ is isomorphic to $Q_{n-1}$, and such that the
restriction of $Q_{0}$ to $\cM_{0}$ has $\spec k$ as its total space. 
We can define a functor from
$R$-algebras to sets that sends each $R$-algebra $A$ to the set of
isomorphism classes of $G$-torsors on the stack $\cM_{A}$. This
functor is easily checked to be limit-preserving (for example, by
using a presentation of $\cM$ and descent for $G$-torsors). So we can
apply Artin's approximation theorem, and conclude that there exists a
$G$-torsor on $\cM$, whose restriction to $\cM_{0}$ has $\spec k$
as its total space.

The total space $\cP$ is an algebraic stack with finite inertia;
furthermore, the inverse image of $\cM_{0}$ in $\cP$ is isomorphic to
$\spec k$. The locus where the inertia stack $\cI_{\cP} \arr \cP$ has
fiber of length larger than $1$ is a closed substack of $\cP$, whose
image in $M = \spec R$ is a closed subscheme that does not contain $p$; hence this locus is empty. So
$\cP$ is an algebraic space (in fact an affine scheme); and this concludes the
proof of the first case. 

{\sc Obtaining a flat morphism and proof of part~\refpart{thm:main}{3}.\enspace} Now we prove a weaker version of the Proposition, with the same statement, except that the morphism $U \arr M$ is only supposed to be flat and finitely presented, instead of \'etale. This is sufficient to prove that part~\refpart{thm:main}{3} holds.

By passing to the algebraic closure of $k$ we may assume that $k$ is algebraically closed.

We claim that there exists a finite extension $k'$ of the residue
field $k(p)$ contained in $k$, such that the object $\xi$ is defined
over $k'$. In fact, it follows from the definition of moduli space
that there exists an object $\eta$ of $\spec \overline{k(p)}$ whose
pullback to $\spec k$ is isomorphic to $\xi$. This $\eta$ gives an
object of the algebraic stack $\spec k(p)\times_{M}\cM$ over $\spec
\overline{k(p)}$, which is finitely presented over $\spec k(p)$, and
any such object is defined over a finite extension $k'$ of $k(p)$. Hence we
may assume that $k$ is a finite extension of $k(p)$. 

There is a flat morphism of finite presentation $M' \arr M$, with a point $q
\in U$ mapping to $p$, such that $k(q) = k$ \cite[$0_{III}$.10.3]{EGA}; hence, by applying the first step to $M' \times_{M} \cM$, there is an \'etale morphism $U \arr M'$ containing $q$ in its image, such that $U\times_{M}\cM$ has the required quotient form.

{\sc The conclusion.\enspace} 
The argument of the proof of the previous case shows that to conclude we only need the following fact.

\begin{proposition}\label{prop:lift-points}
Let $\cM \arr S$ be a tame stack with moduli space $\rho\colon \cM \arr M$, $k$ a field. Given a morphism $\spec k \arr M$, there exists a finite separable extension $k \subseteq k'$ and a lifting $\spec k' \arr \cM$ of the composite $\spec k' \arr \spec k \arr M$.
\end{proposition}

\begin{proof}
We are going to need some preliminaries.

Suppose that $\cX$ and $\cY$ are algebraic stacks over a scheme $S$; consider the stack $\underhom_{S}(\cX, \cY)$ whose sections over an $S$-scheme $T$ are morphisms of\/ $T$-stacks $\cX_{T} \arr \cY_{T}$ (see \cite{homstack}). We will denote by $\homrep_{S}(\cX, \cY)$ the substack whose sections are representable morphisms $\cX_{T} \arr \cY_{T}$.

\begin{lemma}\label{lem:finite-type}
Let $G$ and $H$ be finite and finitely presented flat group schemes over a locally noetherian scheme $S$. Then the stacks $\underhom_{S}(\cB_{S}G, \cB_{S}H)$ and $\homrep_{S}(\cB_{S}G, \cB_{S}H)$ are finitely presented over $S$.
\end{lemma}

\begin{proof}
By standard arguments, we may assume that $S$ is affine and finitely generated over $\ZZ$.

Consider the fppf sheaf $\underhom_{S}(G,H)$ whose sections over an $S$-scheme $T$ are homomorphism of group schemes $G_{T} \arr H_{T}$. This is a locally closed subscheme of the Hilbert scheme of $G\times_{S}H$ over $S$, hence it is of finite type. Furthermore, consider the subsheaf $\hominj_{S}(G, H)$ consisting of universally injective homomorphisms $G_{T} \arr H_{T}$ (a homomorphism is universally injective if it is injective as a homomorphism of sheaves in the fppf topology, and it stays such after arbitrary base change on $T$). We claim that the embedding $\hominj_{S}(G, H) \subseteq \underhom_{S}(G, H)$ is represented by an open subscheme of $\underhom_{S}(G, H)$.

To see this, consider the universal object $\phi\colon G_{T} \arr H_{T}$, where we have set $T \eqdef \underhom_{S}(G, H)$, and its kernel $K \subseteq G_{T}$. The sheaf $K$ is a (not necessarily flat) finite group scheme on $T$. Formation of $K$ commutes with base change on $T$; hence the fiber product $\hominj_{S}(G, H)$ is represented by the open subscheme of $T$ consisting of point $t \in T$ such that $K_{t}$ has length~$1$.

There is a natural action by conjugation of $H$ on $\underhom_{S}(G, H)$, preserving $\hominj_{S}(G, H)$. We claim that $\underhom_{S}(\cB_{S}G, \cB_{S}H)$ is isomorphic to the quotient stack $[\,\underhom_{S}(G, H)/H]$, and that $\homrep_{S}(\cB_{S}G, \cB_{S}H)$ is the open substack $[\,\hominj_{S}(G, H)/H]$. This is obviously enough to prove the statement.

Let us start by producing a functor
   \[
   \Phi\colon \underhom_{S}(\cB_{S}G, \cB_{S}H) \arr [\,\underhom_{S}(G, H)/H].
   \]
For each $S$-scheme $T$ and each base-preserving functor $F\colon \cB_{T}G_{T} \arr \cB_{T}H_{T}$ we need to exhibit an $H$-torsor $Q \arr T$ and an $H$-equivariant morphism $Q \arr \underhom_{S}(G, H)$. The torsor $Q \arr T$ is simply the image via $F$ of the trivial $G$-torsor $G_{T} \arr T$. Furthermore, the functor $F$ induces a homomorphism from $\underaut_{S}(G_{Q} \arr Q) = G_{Q}$ to the automorphism group scheme of the image of $G_{Q} \arr Q$ in $\cB_{Q}H_{Q}$. Since this image is the pullback of $Q$ to $Q$, which is canonically a trivial torsor, its automorphism group scheme is $H_{Q}$. This defines an object of $\underhom_{S}(G,H)(Q)$, hence a morphism $Q \arr \underhom_{S}(G,H)$; this is easily seen to be $H$-equivariant. This gives an object of $[\,\underhom_{S}(G, H)/H](Q)$: the functor $\Phi$ will send $F$ to this object. We leave it to the reader to complete the definition of $\Phi$ by defining its action on arrows.

The functor
   \[
   \Psi\colon [\,\underhom_{S}(G, H)/H] \arr \underhom_{S}(\cB_{S}G, \cB_{S}H)
   \]
is defined as follows. Take an object of $\underhom_{S}(G, H)/H](T)$, that is, a principal $H$-bundle $Q \arr T$, with an $H$-equivariant morphism $\phi\colon Q \arr \underhom_{S}(G, H)$. We need to define a base-preserving functor $\cB_{T}G_{T} \arr \cB_{T}H_{T}$. Suppose that $T'$ is a $T$-scheme and $P \arr T'$ is a $G$-torsor. We can define an action of $G$ on $P \times_{T}Q$ by the formula
   \[
   (p, q)g = \bigl(pg, q\phi(q)(g)\bigr).
   \]
Using the formula
   \[
   \phi(qh) = h^{-1}\phi(q)h
   \]
which expresses the fact that $\phi$ is $H$-equivariant, we check that this action of $G$ commutes with the obvious action of $H$ defined by $(p,q)h = (p, qh)$. The resulting action of $G\times H$ is free and $(P \times_{T}Q)/(G \times_{S}H) = T'$, so $(P \times_{T}Q)/G \arr P/G = T'$ is an $H$-torsor. 

We set $\Psi(Q, \phi)(P \arr T') = (P \times_{T}Q)/G \arr T'$. This defines the action of $\Psi$ on objects; the definition of the action of $\Psi$ on arrows is straightforward.

Let us verify that $\Phi\Psi$ is isomorphic to $\id_{[\,\underhom_{S}(G, H)/H]}$. Let $T$ be an $S$-scheme, $Q\arr T$ an $H$-torsor, $\phi\colon Q \arr \underhom_{S}(G, H)$ an $H$-equivariant maps. Set $\Phi\Psi(Q,\phi) = (Q', \phi')$; we need to produce an isomorphism between $(Q, \phi)$ and $(Q', \phi')$. 

By definition, $Q'$ is the quotient $(G\times_{S}Q)/G$, where the action of $G$ on $G\times_{S}Q$ is defined by the formula $(g, q)g_{1} = \bigl(gg_{1}, q\phi(q)(g_{1})\bigr)$. It is immediate to verify that the morphism $G\times_{S}Q \arr Q$ is $G$-invariant and $H$-equivariant; hence the induced morphism $Q' = (G\times_{S}Q)/G \arr Q$ is an isomorphism of $H$-torsors. We leave it to the reader to check that the homomorphism $\phi'\colon Q' \arr \underhom_{S}(G, H)$ induced by $\Psi(Q,\phi)$ corresponds to $\phi$; this gives the natural isomorphism $(Q, \phi) \simeq (Q', \phi')$.

Next we need to check that $\Psi\Phi$ is isomorphic to $\id_{\underhom_{S}(\cB_{S}G, \cB_{S}H)}$. Let $F\colon \cB_{T}G \arr \cB_{T}H$ be a base-preserving functor. Set $\Phi F = (Q \to T, \phi)$; by definition, $Q \arr T = F(G_{T} \arr T)$. We need to produce a canonical isomorphism between $F$ and $\Psi(Q \to T, \phi)$. Let $T'$ be a $T$-scheme, $P \arr T'$ a $G$-torsor; we need an isomorphism of $H$-torsors between $Q' \arr T \eqdef F(P \to T')$ and $\Psi(Q \to T)(P \to T')$. By definition, $\Psi(Q \to T)(P \to T')$ is the $H$-torsor $(P\times_{T}Q)/G \arr T'$. There is a morphism $P\times_{T}Q \arr Q'$ which is defined as follows. If $T''$ is a $T'$-scheme, a morphism $T'' \arr P$ corresponds to a morphism $G_{T''} \arr P$ of $G$-torsors over $T''$. This in turn induces a morphism of $H$-torsors
   \[
   T''\times_{T} Q = F(G_{T''} \to T'') \arr T''\times_{T'}Q'.
   \]
This construction defines a morphism $P\times_{T}Q \arr Q'$. This is easily checked to be $G$-invariant and $H$-equivariant: hence it defines the desired isomorphism of $H$-torsors $(P\times_{T}Q)/G \simeq Q'$.

We have left to prove that $\homrep_{S}(\cB_{S}G, \cB_{S}H)$ is isomorphic to the quotient $[\,\hominj_{S}(G, H)/H]$, or, equivalently, that the inverse image of the substack $\homrep_{S}(\cB_{S}G, \cB_{S}H) \subseteq \underhom_{S}(\cB_{S}G, \cB_{S}H)$ into $\underhom_{S}(G, H)$ equals $\hominj_{S}(G, H)$. But this follows immediately from the well known fact that a morphism $F\colon \cB_{T}G_{T} \arr \cB_{T}H_{T}$ is representable if and only if for any $T' \arr T$ and any object $\xi$ of $\cB_{S}G(T')$, the induced homomorphism from $\underaut_{T'}(\xi) \arr \underaut_{T'}\bigl(F(\xi)\bigr)$ is injective.
\end{proof}

Let $k$ be a field, $R$ be an artinian local $k$-algebra with residue field $k$, $G$ a \lr group scheme acting on $R$. Set $\cM = [\spec R/G]$, and assume that the moduli space of $\cM$ is $\spec k$ (this is equivalent to assuming that $R^{G} = k$). We have a natural embedding $\cB_{k}G = [\spec k/G] \subseteq [\spec R/G] = \cM$.

\begin{lemma}\label{lem:into-reduced}
If $T$ is a $k$-scheme, any representable morphism of $k$-stacks $\rho :\cB_{k}G\times_{\spec k}T \arr \cM$ factors through $\cB_{k}G \subseteq \cM$.
\end{lemma}

\begin{proof}
Let $P$ be the pullback of $\spec R \arr \cM$ to $\cB_{k}G\times_{\spec k}T$; then $P$ is an algebraic space (since $\rho $ is representable) with an action of $G$, such that the morphism $P \arr \spec R$ is $G$-equivariant. We claim that the composite $P \arr \cB_{k}G \times T \arr T$ is an isomorphism. 

Since it is finite and flat it enough to prove that is an isomorphism when pulled back to a geometric point $\spec \Omega \arr T$, were $\Omega$ is an algebraically closed field; so we may assume that $T = \spec \Omega$. Choose a section $\spec \Omega \arr P$: since there is a unique morphism $\spec \Omega \arr \cB_{\Omega}G_{\Omega}$ over $\Omega$, we get a commutative diagram
   \[
   \xymatrix{
   \spec\Omega \ar[rr]\ar[rd]&& P\ar[ld]\\
   &\cB_{\Omega}G_{\Omega}&
   }
   \]
Since both $\spec \Omega \arr \cB_{\Omega}G_{\Omega}$ and $P \arr \cB_{\Omega}G_{\Omega}$ are $G$-torsors, the degrees of both over $\cB_{\Omega}G_{\Omega}$ equal the order of $G$; hence $\spec \Omega \arr P$ is an isomorphism, and $P \arr \spec \Omega$ is its inverse.

Thus, since the composite $P \arr \cB_{k}G\times_{\spec k}T \arr T$ is $G$-equivariant, this means that the action of $G$ on $P$ is trivial. The morphism $P \arr \spec R$ corresponds to a ring homomorphism $R \arr \cO(P)$, which is $G$-equivariant, and the action of $G$ on $\cO(P)$ is trivial. But if $\frm$ is the maximal ideal of $R$, there is a splitting of $G$-modules $R \simeq \frm \oplus k$; and $\frm^{G} = 0$, because $R^{G} = k$. So $\frm$ is a sum of non-trivial irreducible representations, since $G$ is \lr, and any $G$-equivariant linear map $\frm \arr \cO(P)$ is trivial. So $P \arr \spec R$ factors through $\spec k$, so $\rho :\cB_{k}G\times_{\spec k}T \arr \cM$ factors through $[\spec k/G] = \cB_{k}G$, as claimed.
\end{proof}

Let us prove Proposition~\ref{prop:lift-points}. Since $\cM$ is limit-preserving, it is sufficient to show that any morphism $\spec k \arr M$, where $k$ is a separably closed field, lifts to $\spec k \arr \cM$. 

Notice that Corollary~\refall{2.3}{2} can now be applied, since its proof only requires part~\refpart{thm:main}{3} of the theorem, which we have just verified. Hence the moduli space of $\spec k \times_{M}\cM$ is $\spec k$, and we can base change to $\spec k$, and assume that $M = \spec k$.

Let $k \subseteq k'$ be a finite field extension such that $\cM(k')$ is non-empty. Pick an object $\xi \in \cM(k')$, and set $G_{k'} = \underaut_{k'}(\xi)$. After extending $k'$, we may assume that $G_{k'}$ is of the form $H_{k'} \ltimes \Delta_{k'}$, where $\Delta_{k'}$ is a diagonalizable group scheme whose order is a power of the characteristic of $k$ and $H_{k'}$ is a constant tame group scheme. There exist unique group schemes $\Delta$ and $H$, respectively diagonalizable and constant, whose pullbacks to $\spec k'$ coincide with $\Delta_{k'}$ and $H_{k'}$; furthermore, the action of $H_{k'}$ on $\Delta_{k'}$ comes from a unique action of $H$ on $\Delta$. We set $G = H \ltimes \Delta$: this $G$ is a group scheme on $\spec k$ inducing $G_{k'}$ by base change.

\begin{lemma}
The stack $\homrep_{k}(\cB_{k}G, \cM)$ is of finite type over $k$.
\end{lemma}

\begin{proof}
It is enough to prove the result after base changing to $k'$; we can therefore assume that $\cM = [\spec R/G]$, where $R$ is an artinian $k$-algebra with residue field $k$, because of the first part of the proof. Then by Lemma~\ref{lem:into-reduced} the stack $\homrep_{k}(\cB_{k}G, \cM)$ is isomorphic to $$\homrep_{k}(\cB_{k}G, \cB_{k}G),$$ which is of finite type by Lemma~\ref{lem:finite-type}.
\end{proof}

The morphism $t\colon \spec k \arr \cB_{k}G$ corresponding to the trivial torsor induces a morphism
   \[
   F\colon \homrep_{k}(\cB_{k}G, \cM) \arr \cM
   \]
by composition with $t$.

Consider the scheme-theoretic image $\underline{\cM} \subseteq \cM$ of the morphism $F$: this is the smallest closed substack of $\cM$ with the property that $F^{-1}(\underline{\cM}) = \homrep_{k}(\cB_{k}G, \cM)$. Its sheaf of ideals is the kernel of the homomorphism $\cO_{\cM} \arr F_{*}\cO_{\homrep_{k}(\cB_{k}G, \cM)}$.

\begin{lemma}
We have
   \[
   \spec k' \times_{\spec k} \underline{\cM} = \cB_{k}G_{k'}
   \subseteq \cM_{k'},
   \]
where $\cB_{k'}G_{k'}$ is embedded in $\cM_{k'}$ as the residual gerbe of $\xi$. 
\end{lemma}

\begin{proof}
By the first part of the proof, we can write $\cM_{k'}$ in the form $[\spec R/G_{k'}]$, where $R$ is an artinian $k$-algebra with residue field $k'$. Formation of scheme-theoretic images commutes with flat base change, hence we need to show that the scheme-theoretic image of the morphism
   \[
   F_{k'}\colon \homrep_{k'}(\cB_{k'}G_{k'}, \cM_{k'}) \arr \cM_{k'}
   \]
is equal to $\cB_{k'}G_{k'}$; or, equivalently, that for any morphism
   \[
   g\colon T \arr \homrep_{k'}(\cB_{k'}G_{k'}, \cM_{k'})
   \]
the composite $F_{k'} \circ g\colon T \arr \cM_{k'}$ factors through $[\spec k'/G_{k'}]$. This follows from Lemma~\ref{lem:into-reduced}.
\end{proof}

Now we can replace $\cM$ with $\underline{\cM}$, and assume that $\cM_{k'}$ is $\cB_{k'}G_{k'}$. It follows that $\cM$ is a gerbe in the fppf topology over $\spec k$.

Next we define an \'etale gerbe $\cN$, with a morphism $\cG \arr \cN$.

For any $k$-scheme $T$ and any object $\xi \in \cM(T)$, the automorphism group scheme $G_{\xi} \arr T$ is \lr; let
   \[
   1 \arr \Delta_{\xi} \arr G_{\xi} \arr H_{\xi} \arr 1
   \]
be the connected \'etale sequence of $G_{\xi}$. More concretely, $\Delta_{\xi}$ is the subfunctor of $G_{\xi}$ of automorphisms whose order is a power of the characteristic of $k$. If $f\colon T_1 \arr T$ is a morphism of schemes, then $G_{f^{*}\xi} = T_1 \times_{T} G_{\xi}$ (this is a general property of fibered categories), and $\Delta_{f^{*}\xi} = T_1 \times_{T} \Delta_{\xi}$. We define $\cN$ to be the stack $\cM\thickslash \Delta$, whose existence is assured by Theorem~\ref{thm:rigidification}. We claim that $\cN$ is a Deligne--Mumford stack. It is enough to check that for any algebraically closed field $k$ and any object $\overline{\xi}$ of $\cN(\spec k)$, the automorphism group scheme $\underaut_{k}(\overline{\xi})$ is reduced. However, since the morphism $\cM \arr \cN$ is of finite type, the object $\overline{\xi}$ comes from an object $\xi$ of $\cM(\spec k)$, and we have
   \[
   \underaut_{k}(\overline{\xi}) = G_{\xi}/\Delta_{\xi};
   \]
this is reduced by definition.

The pullback of $\cN$ to $\spec{k'}$ is $\cB_{k'}(G_{k'}/\Delta_{k'})$ (recall that $\Delta_{k'}$ is the connected component of the identity in $G_{k'}$); so $\cN$ is an \'etale gerbe over $\spec k$. Since $k$ is separably closed, there is a $k$-morphism $\spec k \arr \cN$. We can replace $\cM$ by $\cM \times_{\cN} \spec k$, so that $\cM_{k'} = \cB_{k'}\Delta_{k'}$.

In this case, we claim that $\cM$ is banded by the diagonalizable group $\Delta \arr \spec k$ (recall that we have defined this as the diagonalizable group scheme whose pullback to $\spec k'$ is $\Delta_{k'}$). In fact, since $\cM$ is a gerbe, and all of its objects have abelian automorphism groups, then the automorphism group schemes descend to a group scheme over $\spec k$, whose pullback to $\spec k'$ is $\Delta_{k'}$. So this group scheme is a form of $\Delta$ in the fppf topology; but the automorphism group scheme of $\Delta$ is constant, so this form is in fact a form in the \'etale topology, and it is trivial.

The class of the gerbe $\cM$ banded by $\Delta$ is classified by the cohomology group $\H^{2}\fppf(\spec k, \Delta)$.

\begin{lemma}
If $\Delta$ is a diagonalizable group scheme over a separably closed field $k$, we have $\H^{2}\fppf(\spec k, \Delta) = 0$.
\end{lemma}

\begin{proof}
The group scheme $\Delta$ is a product of groups of the form $\mmu_{n}$, so it is enough to consider the case $\Delta = \mmu_{n}$. Then the result follows from the Kummer exact sequence of fppf sheaves
   \[
   0 \arr \mmu_{n} \arr \GG_{\rmm}
   \stackrel{\times n}\longrightarrow \GG_{\rmm} \arr 0
   \]
and the fact that $\H^{i}\fppf(\spec k, \GG_{\rmm}) = \H^{i}\et(\spec k, \GG_{\rmm}) = 0$ for $i>0$ (see \cite[III, Theorem~3.9]{milne}).
\end{proof}
\noqed\end{proof}
\noqed\end{proof}

This concludes the proofs of Propositions \ref{prop:lift-points} and \ref{prop:stronger-version} and of Theorem~\ref{thm:main}.\end{proof}

\appendix
\section{Rigidification}\label{sec:rigidification}

In this section we discuss the notion of \emph{rigidification,} where a subgroup $G$ of inertia is ``removed". This was studied in \cite{A-C-V,Romagny,B-N} when $G$ is in the center of inertia, the general case briefly mentioned in \cite{Kresch}.

In what follows, when we refer to a sheaf on a scheme $T$ this will be a sheaf on the fppf site of $T$.

Let $S$ be a scheme, or an algebraic space, and let $\cX \arr S$ be a locally finitely presented algebraic stack. We will denote by $\cI\cX \arr \cX$ the inertia stack.

Suppose that $G \subseteq \cI\cX$ is a flat finitely presented subgroup stack. This means that it is closed substack, the neutral element section $\cX \arr \cI\cX$ factors through $G$, the inverse morphism $\cI\cX \arr \cI\cX$ carries $G$ into itself, and the composition morphism $\cI\cX\times_{\cX}\cI\cX \arr \cI\cX$ carries $G\times_{\cX}G$ into $G$. Let $\xi$ be an object of $\cX$ over some $S$-scheme $T$; this corresponds to a morphism $T \arr \cX$. The pullback $T\times_{\cX}\cI\cX$ is, canonically, the group scheme $\underaut_{T}(\xi)$; the pullback of $G$ to $T$ gives a flat group subscheme $G_{\xi} \subseteq \underaut_{T}(\xi)$. If $\xi \arr \xi'$ is an arrow in $\cX$ mapping a morphism of $S$-schemes $T \arr T'$, we have a canonical isomorphism $\underaut_{T'}(\xi') \simeq T'\times_{T}\underaut_{T}(\xi)$; this restricts to an isomorphism $G_{\xi'} \simeq T'\times_{T} G_{\xi}$.

Conversely, assume that for each object $\xi$ of $\cX$ we have a subgroup scheme $G_{\xi} \subseteq \underaut_{T}(\xi)$ which is flat and finitely presented over $T$, such that for each arrow $\xi \arr \xi'$ mapping to a morphism $T \arr T'$ the canonical isomorphism $\underaut_{T'}(\xi') \simeq T'\times_{T}\underaut_{T}(\xi)$ carries $G_{\xi'}$ isomorphically onto $T'\times_{T} G_{\xi}$. Then there is a unique flat finitely presented subgroup stack $G \subseteq \cI\cX$ such that for any object $\xi$ of $\cX(T)$ the pullback of $G$ to $T$ coincides with $G_{\xi}$.

Notice that this condition implies that each $G_{\xi}$ is normal in $\underaut_{T}(\xi)$. In fact, if $\xi$ is an object of $\cX(T)$ and $u$ is in $\underaut_{T}(\xi)$, the automorphism of  $\underaut_{T}(\xi)$ induced by the arrow $u\colon \xi \arr \xi$ is conjugation by $u$, and carries $G_{\xi}$ into itself.

\begin{theorem}\label{thm:rigidification}
There exists a locally finitely presented algebraic stack $\cX\thickslash G$ over $S$ with a morphism $\rho\colon \cX \arr \cX\thickslash G$ with the following properties.

\begin{enumeratea}

\item $\cX$ is an fppf gerbe over $\cX\thickslash G$.

\item For each object $\xi$ of $\cX(T)$, the homomorphism of group schemes
   \[
   \underaut_{T}(\xi) \arr \underaut_{T}\bigl(\rho(\xi)\bigr)
   \]
is surjective with kernel $G_{\xi}$.

\end{enumeratea}

Furthermore, if $G$ is finite over $\cX$ then $\rho$ is proper; while if $G$ is \'etale then $\rho$ is also \'etale.
\end{theorem}

\begin{remark}
It is not hard to show that these properties characterize $\cX\thickslash G$ uniquely up to equivalence.

Furthermore, the sheaf of groups $G$ on $\cX$ descends to a band $\underline{G}$ over $\cX\thickslash G$, and $\cX$ is banded by $\underline{G}$ over $\cX\thickslash G$. Using this one can give a modular interpretation of $\cX\thickslash G$ in a spirit similar to that of \cite[Proposition~C.2.1]{A-G-V}.
\end{remark}

\begin{example}
Suppose that $S$ is a locally noetherian scheme, $\cX \arr S$ a regular Deligne--Mumford stack with finite inertia, locally of finite type over $S$. We can take as $G$ the closure of the fibers of $\cI\cX$ over the generic points of the irreducible components of the moduli space of $\cX$. The morphism $G \arr \cX$ is easily seen to be \'etale; hence we can construct a rigidification $\cX\thickslash G$, which is a regular Deligne--Mumford stack with trivial generic stabilizers. Therefore every regular Deligne--Mumford is an \'etale gerbe over a regular Deligne--Mumford stack with trivial generic stabilizers.

This fact was used, for example, in \cite{kv}, without adequate justification.
\end{example}

\begin{proof}
Let us define a fibered category $\overline{\cX}$ over $S$ in the following fashion. 

The objects of $\overline{\cX}$ are the objects of $\cX$.

Let $f\colon T \arr T'$ be a morphism of $S$-schemes, $\xi$ and $\xi'$ be objects of $\cX$ mapping to $T$ and $T'$ respectively. The set of arrows $\xi \arr \xi'$ mapping to $f$ is in a canonical bijective correspondence with the global sections of the sheaf $\underisom_{T}(\xi, f^{*}\xi')$. There is a faithful right action of $G_{\xi}$ on the algebraic space $\underisom_{T}(\xi, f^{*}\xi')$ 

Notice the following fact. There is also a left action of $T\times_{T'}G_{\xi'}$ on $\underisom_{T}(\xi, f^{*}\xi')$ by composition. This commutes with the action of $G_{\xi}$. Furthermore, it is easy to see that the induced action of $T\times_{T'}G_{\xi'}$ on the quotient $\underisom_{T}(\xi, f^{*}\xi')/G_{\xi}$ is trivial; hence the quotient $\underisom_{T}(\xi, f^{*}\xi')/G_{\xi}$ coincides with the double quotient $G_{\xi'}\backslash\underisom_{T}(\xi, f^{*}\xi')/G_{\xi}$.

	We define an arrow $\xi \arr \xi'$ in $\overline{\cX}$ mapping to $f$ as the set of global sections of the quotient sheaf $\underisom_{T}(\xi, f^{*}\xi')/G_{\xi}$.

Let $g\colon T' \arr T''$ be a morphism of $S$-schemes, and let $\xi''$ be an object of $\cX$ mapping to $T''$. Composition defines a morphism of algebraic spaces
   \[
   \underisom_{T}(\xi, f^{*}\xi')\times_{T'} \underisom_{T'}(\xi', g^{*}\xi'')
   \arr \underisom_{T}(\xi, (gf)^{*}\xi'')
   \]
which descends to a morphism
   \[
   \underisom_{T}(\xi, f^{*}\xi')/G_{\xi}
      \times_{T'} \underisom_{T'}(\xi', g^{*}\xi'')/G_{\xi'}
   \arr \underisom_{T}(\xi, (gf)^{*}\xi'')/G_{\xi};
   \]
this defines the composition in $\overline{\cX}$.

The fibered category $\overline{\cX}$ is a prestack, but not a stack. We define $\cX\thickslash G$ to be the stack over $S$ associated with $\overline{\cX}$. The morphism $\rho\colon \cX \arr \cX\thickslash G$ is the composite of the two obvious morphisms $\cX \arr \overline{\cX}$ and $\overline{\cX} \arr \cX\thickslash G$.

Let $T$ be and $S$-scheme, $\overline{\xi}$ and $\overline{\eta}$ two objects in $(\cX\thickslash G)(T)$. We claim that the sheaf $\isom_{T}(\overline{\xi}, \overline{\eta})$ is an algebraic space. This is a local problem in the fppf topology, by a result of Artin (\cite[Corollaire~10.4.1]{L-MB}), so we may assume that $\overline{\xi}$ and $\overline{\eta}$ are the images of two objects $\xi$ and $\eta$ of $\cX(T)$. In this case $\isom_{T}(\overline{\xi}, \overline{\eta})$ coincides with the quotient $\underisom_{T}(\xi, f^{*}\xi')/G_{\xi}$, which is an algebraic space, again by Artin's result (\cite[Corollaire~10.4]{L-MB}).

Now we check that the morphism $\rho\colon \cX \arr \cX\thickslash G$ is represented by smooth algebraic stacks which are fppf gerbe. More precisely, we now show that, given a morphism $T \arr \cX\thickslash G$, where $T$ is a scheme, the pullback $T \times_{\cX\thickslash G} \cX$ is fppf locally isomorphic to a stack of the form $\cB_{T}G_{\xi}$. This implies the statement of the theorem. If $U \arr \cX$ is a smooth morphism, then the composition $U \arr \cX\thickslash G$ is also smooth, which shows that $\cX\thickslash G$ is an algebraic stack. Furthermore $\cB_{T}G_{\xi}$ is proper over $T$ when $G_{\xi}$ is finite, and \'etale when $G_{\xi}$ is \'etale.

We may assume that the object $\overline{\xi}$ of $(\cX\thickslash G)(T)$ corresponding to the morphism $T \arr \cX\thickslash G$ comes from an object $\xi$ of $\cX(T)$. An object of the fiber product $T \times_{\cX\thickslash G} \cX$ consist of the following data:

\begin{enumeratea}

\item A morphism of schemes $f\colon T' \arr T$;

\item An object $\eta$ of $\cX(T')$;

\item A section
   \[
   T' \arr \underisom_{T'}\bigl(\rho(h), \overline{\xi}\bigr)
      = \underisom_{T'}(\eta,\xi)/G_{\xi}.
   \]

\end{enumeratea}

From this data we obtain a $G_{\xi}$-torsor $P \arr T'$, where
   \[
   P \eqdef
      T' \times_{\underisom_{T'}(\eta,\xi)/G_{\xi}}\underisom_{T'}(\eta,\xi).
   \]
By examining what happens to arrows it is easy to show that this construction extends naturally to a functor $T \times_{\cX\thickslash G} \cX \arr \cB_{T}G_{\xi}$. We will show that this is an equivalence. The inverse $\cB_{T}G_{\xi} \arr T \arr \cX\thickslash G$ is defined as follows. Let $T'$ be a scheme over $T$ and $\phi\colon P \arr T'$ be a $G_{\xi}$-torsor. The action of $G_{\xi}$ on $\xi$ gives to an action of $G_{\xi}$ on the pullback $\xi_{T}$; by \cite[Theorem~4.46]{descent} we obtain an object $\eta$ of $\cX(T')$. The functor $\cB_{T}G_{\xi} \arr T \arr \cX\thickslash G$ sends $P \arr T'$ into $\eta$.
We leave it to the reader to check that this gives an inverse to the functor above. This concludes the proof.

\end{proof}

\bibliographystyle{amsplain}
\bibliography{mrabbrev,VistoliRefs}

\end{document}